\documentstyle[12pt,amstex,amssymb,amsthm]{article}

\newcommand{\C}{\mbox{\rm \,l\kern-0.52em C}}
\newcommand{\Ce}{\rm \,l\kern-0.35em C}
\newcommand{\R}{{\rm l}\!{\rm R}}
\newcommand{\Z}{{\sf Z}\!\!{\sf Z}}

\newcommand{\N}{{\rm l}\!{\rm N}}

\newtheorem{theorem}{Theorem}[section]
\newtheorem{deflem}[theorem]{Definition and Lemma}
\newtheorem{prop}[theorem]{Proposition}

\newtheorem{problem}[theorem]{Problem}

\newtheorem{lemma}[theorem]{Lemma}
\newtheorem{remark}[theorem]{Remark}

\renewenvironment{proof}{{\bf Proof:}}{\mbox{}\hfill $\Box$}
\theoremstyle{definition}
\newtheorem{definition}[theorem]{Definition}
\title{Characters of Fredholm modules and a problem of Connes}
\author{Michael Puschnigg}
\date{}

\begin{document}
\maketitle

A large part of research in noncommutative geometry is devoted to various far
reaching generalizations of the Atiyah-Singer index theory of elliptic
operators on compact manifolds.  Following Kasparov \cite{Ka}, the notion of
an elliptic operator on a manifold is replaced by that of a Fredholm module (or
$K$-cycle) over an algebra of operators on Hilbert space. A number of index
theorems for Fredholm modules associated to geometric data have been obtained
by Kasparov \cite{Ka}, Connes (see \cite{Co2} and the literature cited
therein) and others. In order to handle concrete index problems it is not
only necessary to dispose of an index theorem but also to provide an index
formula which allows the explicit calculation of indices. The classical index formula 
of Atiyah and Singer is obtained from the index
theorem by applying the Chern character in topological $K$-homology. This
motivates the search for character formulas of Fredholm modules ($K$-cycles)
that define a Chern character on the $K$-homology groups of an algebra of
operators.  

It was this search for a Chern
character in $K$-homology which led A.~Connes to the invention of cyclic
cohomology \cite{Co1}. He obtained various explicit character formulas
\cite{Co1}, \cite{Co}, \cite{Co2} which depend on   the degree of analytic
regularity (summability) of the given Fredholm module.  Whereas the classical
index problems of Atiyah-Singer are all finitely summable in the sense of
\cite{Co1} there are many examples of Fredholm modules over (noncommutative)
algebras which are infinite-dimensional (not finitely summable). Characters of
Fredholm modules have been calculated in many finitely summable cases, but as
far as the author knows the character of a  $\theta$-summable, infinite
dimensional (unbounded) Fredholm module \cite{Co}  has not yet been determined
in a single case. 

Typical examples of infinite dimensional (unbounded) Fredholm
modules are modules over dense subalgebras of the reduced group $C^*$-algebra
$C^*_r(\Gamma)$ of a discrete nonamenable group $\Gamma$ \cite{Co5},
\cite{Co2}. A particularly interesting example is presented  by Connes in
\cite{Co3}, \cite{Co2}, where he constructs an infinite dimensional unbounded
Fredholm module ${\cal E}_\Gamma$ over the group ring of a discrete subgroup
$\Gamma$ of a real semisimple Lie group $G$. This module is closely related to
Kasparov's $\gamma$-element \cite{Ka1}. In \cite{Co3}, \cite{Co2}
Connes makes several steps towards the calculation of the character of ${\cal
E}_\Gamma$ and predicts that it should be cohomologous to the canonical
trace on $C^*_r(\Gamma)$. He poses the verification as a problem
\cite{Co3},p.83 and notes that a positive solution would imply the
Kadison-Kaplansky idempotent conjecture for $\Gamma$ \cite{Co1}, \cite{Co2}.

In this paper we solve Connes' problem for cocompact discrete subgroups
$\Gamma$ of semisimple Lie groups $G$ of real rank one. (This restriction is
necessary because we have to make use of Jolissaint's rapid decay property 
for $\Gamma$ \cite{CM}.) To be more precise we prove the 
\begin{theorem}
Let $\Gamma$ be a discrete cocompact subgroup of the semisimple Lie group $G$
of real rank one. Let ${\cal E}_\Gamma$ be Connes' unbounded
$\theta$-summable Fredholm module over ${\cal A}(\Gamma)\subset C^*_r(\Gamma)$
\cite{Co2} which represents the $\gamma$-element
$\gamma_\Gamma \in KK^\Gamma(\C,\C)$ \cite{Ka1}. Then the entire
cyclic character cocycle \cite{Co}  
$$
Ch_{\epsilon}({\cal E}_\Gamma)\,\in\,CC_{\epsilon}^0({\cal A}(\Gamma))
$$
of ${\cal E}_\Gamma$ is cohomologous in local cyclic
cohomology \cite{Pu2} to the canonical trace on ${\cal A}(\Gamma)$
$$
[Ch_{\epsilon}({\cal E}_\Gamma)]\,=\,[\tau_{can}]\,\in\,HC_{loc}^0({\cal
A}(\Gamma)) 
$$
\end{theorem}

There are two main ingredients in the proof of the theorem. The first is a
comparison theorem for character formulas attached to a given
Fredholm module. The second is the partial calculation of the
equivariant bivariant Chern-Connes character of the $\gamma$-element
obtained in \cite{Pu3}. 

The largest part of this paper is concerned with a detailed analysis 
of the relation between Connes' explicit character formula for an unbounded 
$\theta$-summable Fredholm module over an algebra of
operators \cite{Co} and the bivariant Chern-Connes character for bounded 
Kasparov-bimodules over $C^*$-algebras \cite{Pu1}. The
following comparison theorem is the second main result of the paper.

\begin{theorem}
Let ${\cal E}$ be an unbounded $\theta$-summable Fredholm module over the
algebra $A$.  Let $\overline{A}$ be the enveloping $C^*$-algebra of $A$ and
let  $[{\cal E}]\,\in\,KK(\overline{A},\C)$ be the $K$-homology class
corresponding to the given module (it may be represented by any bounded
Fredholm module attached to $\cal E$ by functional calculus). Denote by
$Ch_\epsilon$ Connes'  character of unbounded $\theta$-summable Fredholm
modules in entire cyclic cohomology and let $ch_{biv}$ be the bivariant
Chern-Connes character on $KK$-theory with values in (bivariant) local cyclic
cohomology. Then the images of the characters in the diagram $$
\begin{array}{ccccc}
HC_{\epsilon}^0(A) & \to & HC_{loc}^0(A) & \leftarrow &
HC_{loc}^0(\overline{A}) \\
 & & & & \\
{[}Ch_{\epsilon}({\cal E}){]} & \to & * & \leftarrow & ch_{biv}({[}{\cal E}{]})
\\
\end{array}
$$
coincide.
\end{theorem}

We make a few comments about this result. The advantage
of Connes' character in entire cyclic cohomology lies in the fact that it
provides a completely explicit character formula. However it turns out 
to be very rigid in several ways. An unbounded
Fredholm module over $A$ defines a $K$-homology class over the enveloping
$C^*$-algebra $\overline{A}$. The entire cyclic cohomology groups of $A$ and
$\overline{A}$ will usually differ from each other. So if two $K$-cycles
define the same $K$-homology class but have different domains it may be
impossible to compare their characters because they lie in completely
different cohomology groups. Even if two homotopic $K$-cycles  have the same
domain  $A\subset\overline{A}$, their characters will not necessarily be
cohomologous. An example is given by the pull-back along  two different
evaluation maps $C[0,1]\longrightarrow\C$ of an unbounded Fredholm module with
nonzero index over $\C$. This happens because the entire cyclic theory is
not a (continuous) homotopy functor, as was shown by Khalkhali \cite{Kh}. 

The bivariant Chern-Connes character is of a very different type.
Its existence is established by an
abstract category theoretic argument based on the axiomatic characterization
of bivariant $K$-theory \cite{Hi}. 
Because the definition of the Chern-Connes character involves excision in local
cyclic cohomology, no explicit character formulas exist. This is the main
drawback of this character. On the other hand the Chern-Connes character 
of a $K$-cycle depends by construction only on its $K$-homology class. 
Moreover this character has excellent functorial properties due to its
multiplicativity with respect to the Kasparov-product on bivariant $K$-theory.

Thus the theorem above provides a link between two characters of very
complementary nature. This explains the interest in such a  result.

I heartly thank Prof.~A.~Connes for helpful discussions and suggestions.

\tableofcontents

\section{Characters of Fredholm modules}

We collect some well known material about Fredholm modules, cyclic cohomology 
and Chern characters of Fredholm modules. The only new result is a description
of Connes' Chern character of an unbounded $\theta$-summable fredholm module
\cite{Co} in terms of the Cuntz-Quillen picture of cyclic cohomology. This is
slightly different from Connes' point of view and emphasizes 
the similarities of Connes' character and the
bivariant Chern-Connes character of Cuntz \cite{Cu1}, \cite{Pu1}, \cite{Pu2}.

{\large{\bf Fredholm modules}}

We recall for the convenience of the reader the definition 
of (un)bounded Fredholm modules over an algebra.

\begin{definition}\cite{Ka},\cite{Co1}.
A bounded even Fredholm module over a $C^*$-algebra $A$ is a quadruple
${\cal E}\,=\,({\cal H},\,\rho,\,\epsilon,\,F)$ consisting of 
\begin{itemize}
\item[$\bullet$] a separable Hilbert space $\cal H$.
\item[$\bullet$] a selfadjoint operator $\epsilon\in{\cal L(H)}$ satisfying 
$\epsilon^2\,=\,1$. It defines a $\Z/2\Z$-grading $\cal H=\cal
H^+\oplus\cal H^-$ of $\cal H$ which is given by the $\pm 1$-eigenspaces of
$\epsilon$.   
\item[$\bullet$] a representation
$\rho:A\,\longrightarrow\,{\cal L(H)}$ which is even in the sense that it
commutes with the grading operator $\epsilon$.
\item[$\bullet$] an operator $F\in{\cal L(H)}$ which satisfies
$$
\begin{array}{cc}
\rho(a)(F^2-1)\in{\cal K(H)}, & \rho(a)(F-F^*)\in{\cal K(H)}, \\
 & \\
\rho(a)(\epsilon F\,+\,F\epsilon)\in{\cal K(H)}, & 
[F,\,\rho(a)]\in{\cal K(H)} \\
\end{array}
\eqno(1.1)
$$
for all $a\in A$, where ${\cal K(H)}$ denotes the ideal of compact operators in
${\cal L(H)}$. An odd Fredholm module is defined similarly by  forgetting the
grading of $\cal H$. 
\end{itemize}
\end{definition}

The set of homotopy classes of even (odd) Fredholm modules over a 
$C^*$-algebra $A$ carries a natural abelian group structure \cite{Ka}. 
These are the $K$-homology groups $K^*(A)\,=\,KK_*(A,\C)$ of $A$ \cite{Ka}. 

\begin{definition}\cite{Ka},\cite{Co1}.
An unbounded even Fredholm module over an algebra $A$ is a quadruple
${\cal E}\,=\,({\cal H},\,\rho,\,\epsilon,\,{\cal D})$ consisting of 
\begin{itemize}
\item[$\bullet$] a separable Hilbert space $\cal H$.
\item[$\bullet$] a selfadjoint operator $\epsilon\in{\cal L(H)}$ satisfying 
$\epsilon^2\,=\,1$. It defines a $\Z/2\Z$-grading of $\cal H$ given by the
$\pm 1$-eigenspaces of $\epsilon$. 
\item[$\bullet$] a representation $\rho:A\,\longrightarrow\,{\cal L(H)}$
which is even in the sense that it commutes with
the grading operator.
\item[$\bullet$] an unbounded selfadjoint operator ${\cal D}$ on $\cal H$ which
satisfies 
$$
\begin{array}{ccc}
\epsilon{\cal D}+{\cal D}\epsilon=0, & 
\rho(a)(1+{\cal D}^2)^{-1}\in{\cal K(H)}, & [{\cal D},\,\rho(a)]\in{\cal L(H)} 
\\ \end{array}
\eqno(1.2)
$$
for all $a\in A$. An odd Fredholm module is defined similarly by 
forgetting the grading of $\cal H$. 
\end{itemize}
\end{definition}

If $\cal E$ is an unbounded Fredholm module with domain $A$, then
we suppose that $A$ is complete with respect to the norm
$$
\parallel a\parallel_A\,=\,\parallel\rho(a)\parallel_{{\cal L(H)}}\,+\,
\parallel[{\cal D},\,\rho(a)]\parallel_{{\cal L(H)}}
\eqno(1.3)
$$
Therefore $A$ is a Banach algebra.

Let ${\cal E}\,=\,({\cal H},\,\rho,\,\epsilon,\,{\cal D})$ be an (even) 
unbounded Fredholm module over $A$. Denote by $\overline{A}$ 
the enveloping $C^*$-algebra of $\rho(A)\subset{\cal L(H)}$ and let $f$ be a
continuous function on the real line such that  $\underset{t\to
+\infty}{lim}\,f(t)\,=\,+1,\,\underset{t\to -\infty}{lim}\,f(t)\,=\,-1$. Then 
the quadruple ${\cal E}'\,=\,({\cal H},\,\rho,\,\epsilon,\,f({\cal D}))$
defines a bounded (even) Fredholm module over $\overline{A}$. Its $K$-homology
class  depends only on $\cal E$ and not on the choice of $f$. It is called the
$K$-homology class associated to the given unbounded Fredholm module.

It is possible to impose certain normalization conditions on bounded 
Fredholm modules without changing their $K$-homology class. So one 
can achieve by a suitable homotopy that 
$$
F^2-1\,=\,F-F^*\,=\,
\epsilon F\,+\,F\epsilon\,=\,0
$$
 We will however demand less and will suppose henceforth
that the considered bounded Fredholm modules satisfy apart from (1.1) 
the normalization condition
$$
\begin{array}{ccc}
F^2=1, & F-F^*\in{\cal K(H)}, &
\epsilon F\,+\,F\epsilon\in{\cal K(H)}
\end{array}
\eqno(1.4)
$$
If ${\cal E}\,=\,({\cal H},\,\rho,\,\epsilon,\,{\cal D})$ is an unbounded even Fredholm module then 
$$
{\cal E}(z)\,=\,({\cal H},\,\rho,\,\epsilon,\,F(z)),\,F(z)\,=\,\frac{{\cal
D}+z^{\frac12}\epsilon}{({\cal
D}^2+z)^{\frac12}},\,z\in\C\backslash\R_- \eqno(1.5)
$$
is a holomorphic family of bounded Fredholm modules satisfying (1.4)
and representing the $K$-homology class of $\cal E$.  Here $z\to z^{\frac12}$
denotes the branch of the square-root function on the domain
$\C\backslash\R_-$ which takes the value $+1$ at 1.

The character formulas for Fredholm modules which we are going to study 
often require the involved operators to satisfy certain regularity
conditions. The notion of finite summability, i.e.
$Trace((1+{\cal D}^2)^{-p})<\infty$ for $p>>0$ \cite{Co1} is well known 
but rather restrictive. For example there exists no finitely summable 
unbounded Fredholm module over a dense subalgebra of the reduced group
$C^*$-algebra of a  nonamenable discrete group \cite{Co5}. Most known examples
of Fredholm modules satisfy however the following much weaker
regularity condition.
  
\begin{definition}\cite{Co}
An unbounded Fredholm module ${\cal E}\,=\,({\cal H},\,\rho,\,\epsilon,\,{\cal D})$ 
is \\ $\theta$-summable if
$$
\begin{array}{cc}
Trace(e^{-t{\cal D}^2})\,<\,\infty, & \forall t>0
\end{array}
\eqno(1.6)
$$
\end{definition}

\begin{remark}
A.~Connes works in \cite{Co2} also with Fredholm modules which satisfy 
$Trace(e^{-t{\cal D}^2})\,<\,\infty$ only for $t>>0$. His character formula
makes perfectly sense for such modules. The methods of this paper apply however
only under the condition that the whole heat semigroup  of $\cal D$ is of trace
class. \end{remark} 

A basic invariant of an even
Fredholm module is its index 
$$
Index({\cal E})\,=\,dim \,\,Ker({\cal D}\vert_{{\cal H}_+})\,-\,dim
\,\,Ker({\cal D}\vert_{{\cal H}_-})
\eqno(1.7)
$$
The index of a Fredholm module is homotopy invariant and depends therefore only
on its $K$-homology class.

{\large{\bf Cyclic cohomology theories}}

We recall some well known facts about cyclic complexes  
taken from Connes \cite{Co},\cite{Co1},\cite{Co2}, Cuntz and Quillen
\cite{CQ},\cite{CQ1},\cite{CQ2} and \cite{Pu1},\cite{Pu2}. We only treat
cohomology because the corresponding
homology groups will play no role in this paper. 

Let $A$ be a complete, locally convex algebra over the complex numbers with 
jointly continuous multiplication. The universal complete, locally convex,
differential graded algebra $\Omega A$ over $A$ is given by   
$$
\begin{array}{cccc} 
\Omega A\,=\,\underset{n=0}{\overset{\infty}{\bigoplus}}\,\Omega^n A & \simeq &
A\,\oplus\,\underset{n=1}{\overset{\infty}{\bigoplus}}
\,\widetilde{A}\otimes_\pi A^{\otimes_\pi n} & \\  & & & \\
a^0da^1\ldots da^n & \leftrightarrow &
a^0\otimes a^1\otimes\ldots\otimes a^n, & n\geq 0 \\  
 & & & \\
da^1\ldots da^n & \leftrightarrow &
1\otimes a^1\otimes\ldots\otimes a^n, & n\geq 1 \\  
\end{array} \eqno(1.8)
$$ 
where $\widetilde{A}$ 
is the algebra obtained from $A$ by adjoining a unit and $\otimes_\pi$ denotes
the projective tensor product. The algebra $\Omega A$ is naturally
$\Z/2\Z$-graded by the decomposition into forms of even and odd degree,
respectively
.  
The basic operators on algebraic differential forms are  
the exterior derivative 
$$ 
\begin{array}{cc} 
d:\Omega^nA\longrightarrow\Omega^{n+1}A, &  
d(a^0da^1\ldots da^n)=da^0da^1\ldots da^n\\ 
\end{array} \eqno(1.9)
$$ 
and the  Hochschild boundary operator  
$$ 
\begin{array}{lr} 
b:\,\Omega^nA\longrightarrow\Omega^{n-1}A, & b(\omega da)= 
(-1)^{\vert\omega\vert}[\omega, a]\\ 
\end{array} \eqno(1.10)
$$ 
They satisfy $b^2\,=\,d^2\,=\,0$. 
The Karoubi operator 
$$ 
\begin{array}{lr} 
\kappa:\Omega^nA\longrightarrow\Omega^nA, &  
\kappa(a^0da^1\ldots da^n)=(-1)^{n-1}da^na^0da^1\ldots da^{n-1} 
\end{array} \eqno(1.11)
$$ 
can be expressed in terms of the two basic operators as 
$\kappa=Id-(d\circ b+b\circ d)$ Therefore it commutes with $b$ and $d$. 
The Connes operator  
$$ 
\begin{array}{cc} 
B:\Omega^nA\longrightarrow \Omega^{n+1}A, & B\,=\, 
\sum_{i=0}^{n}\kappa^i\circ d \\ 
\end{array} \eqno(1.12)
$$  
satisfies $B^2=0$ and $bB+Bb=0$ and commutes with $\kappa$ as well. 
The operators $d$, $b$ and $B$ are of degree one with respect to the
$\Z/2\Z$-grading of $\Omega A$.

Let $\Omega A'\,=\,\prod\,\Omega^n A'$ be the dual space of bounded linear
functionals on $\Omega A$. It is a $\Z/2\Z$-graded vector space. 
The operators $d,b$ and $B$ on $\Omega A$ give rise to dual operators
on $\Omega A'$ which will be denoted by the same letters.
The periodic cyclic cochain complex of $A$ is defined as   
$$ 
\begin{array}{ccc} 
CC_{per}^*(A) & = & (\Phi=\prod\varphi_n:\Omega A=\bigoplus
\Omega^nA\to\C,\varphi_n=0,n>>0,\,b+B) \\  
\end{array} \eqno(1.13)
$$ 

It is a $\Z/2\Z$-graded chain complex. Its cohomology 
$HP^*(A)$ is the periodic cyclic cohomology of $A$. 

Suppose that $A$ is a Banach algebra. Then there are various 
cyclic (co)homology theories which take the topology of $A$ more appropriately
into account than the periodic cyclic theory. 
For a linear functional $\Phi=\prod\varphi_n:\Omega A=\bigoplus
\Omega^nA\to\C$ and a subset $S\subset A$ put
$$
\begin{array}{cc}
\parallel\Phi\parallel_S\,=\,
\underset{n\in\N}{Sup}\,\underset{a^0,\ldots,a^n\in S}{Sup}\,
c_n\,\vert\varphi_n(a^0,\ldots,a^n)\vert,\,\,\,c_{2n}=c_{2n+1}=\frac{1}{n!}
\end{array}
\eqno(1.14)
$$
The entire and analytic cyclic cochain complexes of $A$ are defined
as 
$$
\begin{array}{ccc}
 CC_\epsilon^*(A) & = & (\Phi:\Omega
A\to\C,\,\parallel\Phi\parallel_S\,<\infty,\,S\subset A
\,\,\text{bounded},\,b+B)
\end{array}
\eqno(1.15) 
$$ 
and
$$
\begin{array}{ccc}
CC_{anal}^*(A) & = & (\Phi:\Omega
A\to\C,\,\parallel\Phi\parallel_S\,<\infty,\,S\subset A
\,\,\text{compact},\,b+B)
\end{array}
\eqno(1.16) 
$$
respectively.
The cohomology groups of these $\Z/2\Z$-graded complexes are the entire cyclic
cohomology $HC^*_\epsilon(A)$ and the analytic cyclic cohomology groups
$HC^*_{anal}(A)$ of $A$, repectively. 
Finally there is the local cyclic cohomology theory $HC^*_{loc}$ for Banach
algebras. It is invariant under continuous homotopies which makes it
particularly well behaved for $C^*$-algebras. Its definition will not be
needed in the sequel.

The  cyclic (co)homology theories discussed up to now appear as special cases 
of bivariant homology theories
$HC_*^{\beta}(-,-),\,\beta\,=\,per,\,\epsilon,\,anal,\,loc$. These theories 
are contravariant in the first and covariant in the second variable.
The associated (co)homology theories are obtained by fixing one of the
variables
$$
\begin{array}{cc}
HC_*^{\beta}(-)\,=\,HC_*^{\beta}(\C,-), &
HC^*_{\beta}(-)\,=\,HC_*^{\beta}(-,\C) \\
\end{array}
\eqno(1.17)
$$
The bivariant cohomology theories possess an associative and unital composition
product 
$$
\circ:\,HC_*^{\beta}(A,B)\,\otimes\,HC_*^{\beta}(B,C)\,\longrightarrow\,HC_*^{\beta}(A,C)
\eqno(1.18)
$$
and an associative exterior product
$$
\times:\,HC_*^{\beta}(A,B)\,\otimes\,HC_*^{\beta}(C,D)\,\longrightarrow\,
HC_*^{\beta}(A\otimes_\pi C,B\otimes_\pi D)
\eqno(1.19)
$$
All cyclic theories listed up to now satisfy excision: short exact sequences 
of algebras with bounded linear section induce long exact sequences of 
bivariant cyclic groups in any of the two variables. 
There are obvious forgetful transformations 
$$
HC_*^{anal}\,\longrightarrow\,HC_*^\epsilon\,\longrightarrow\,HC_*^{per}
\eqno(1.20)
$$
of homology and dual transformations
$$
HC^*_{per}\,\longrightarrow\,HC^*_\epsilon\,\longrightarrow\,HC^*_{anal}
\eqno(1.21)
$$ 
of cyclic cohomology theories. Moreover there is a natural forgetful
transformation  
$$
HC_*^{anal}(-,-)\,\longrightarrow\,HC_*^{loc}(-,-)
\eqno(1.22)
$$ 
of bifunctors which is mutiplicative with respect to the composition and
exterior products. All transformations are compatible with excision.

The interest in cyclic (co)homology theories originates from the existence
of a natural Chern character 
$$
ch:\,K_*\,\longrightarrow\,HC_*^{loc}\to
HC_*^{anal}\to HC_*^\epsilon\to HC_*^{per} 
\eqno(1.23)
$$
in operator $K$-theory which is compatible (up to universal constants) with
exterior products and excision. The Chern character turns out to be a special
case  of a bivariant Chern-Connes character which will be discussed in 1.7.

{\large{\bf Characters of Fredholm modules}}

The need of character formulas for Fredholm modules was a basic 
motivation for A.~Connes in his search for a noncommutative differential
geometry \cite{Co1}. Character formulas allow the explicit calculation 
of the index of a Fredholm module in exactly the same way as the
Atiyah-Singer index formula leads to the calculation of indices of elliptic
operators on compact manifolds. A character formula associates 
to a Fredholm module $\cal E$ over $A$ a cyclic cocycle $\check{ch}({\cal
E})$ over $A$ such that the index formula  
$$
\langle\check{ch}({\cal E}),ch(e)\rangle\,=\,Index({\cal E}_e)
$$
holds for every $K$-theory class $[e]\in K_0(A)$. Here ${\cal E}_e$ means 
the Fredholm module $\cal E$ twisted by $[e]$ and $ch$ denotes 
the Chern character on $K$-theory with values in cyclic homology \cite{Co1}. 
Contrary to the Chern character in $K$-theory there is no universal character 
formula working simultaneously for all Fredholm modules. All known explicit
character formulas \cite{Co1},\cite{Co},\cite{CM1}, \cite{Ni}
 work only under certain (sometimes rather restrictive) regularity
assumptions on the involved operators. Moreover the cohomology class of the
character cocycles changes definitely under homotopy. Therefore these
characters cannot descend to $K$-homology. There is a character formula in
analytic cyclic cohomology  \cite{Me} which works for all (bounded) modules
without regularity assumptions but involves an infinite number of auxiliary
choices. It descends in fact to $K$-homology but it is not known whether it
can be extended to the bivariant setting. Finally there is a bivariant
multiplicative Chern-Connes character on  Kasparov's bivariant $K$-theory with
values in bivariant local cyclic cohomology \cite{Pu1}. Its existence is a
consequence of the axiomatic characterization of $KK$-theory
\cite{Hi}, but it is difficult  to obtain explicit formulas for this abstract
character because its  construction makes use of excision in cyclic
cohomology. Thus there are various partial solutions of the problem of finding
a universal character formula. The comparison of such partially defined
characters can be quite tedious  as can be seen in this paper. 

We want to present two of the described characters and 
will study their relation in detail. The first one is A.~Connes' character 
of a $\theta$-summable unbounded Fredholm module in entire cyclic
cohomology \cite{Co}. The other one is the abstract Chern-Connes character
on $K$-homology with values in local cyclic cohomology \cite{Pu1}. For
simplicity we will treat only even Fredholm modules. All results of this paper 
hold for odd modules as well.

{\bf Connes' character of a $\theta$-summable unbounded Fredholm module}

\begin{theorem}\cite{Co}
Let ${\cal E}\,=\,({\cal H},\,\rho,\,\epsilon,\,{\cal D})$ be an even unbounded 
$\theta$-summable Fredholm module over $A$. The functionals 
$$
\begin{array}{lr}
Ch_\epsilon({\cal E})\,=\,(\varphi_{2n})_{n\in\N},
& \varphi_{2n}(a^0,\ldots,a^{2n})\,= 
\end{array}
$$
$$
=\,c_n\,\,Trace'\left(\int_{1-i\infty}^{1+i\infty}
F(z^2)\,a^0\,[F(z^2),\,a^1]\ldots[F(z^2),\,a^{2n}]\,e^{z^2}\,dz\right)
\eqno(1.24)
$$
$$
\begin{array}{lcr}
c_n=\frac{1}{\sqrt{\pi}\,i}\,(-\frac{1}{4})^n\,\frac{1}{n!}, &
Trace'(T)=\frac12 Trace(T+\epsilon T\epsilon), & F(z)=\frac{{\cal
D}+\epsilon\,z^{\frac12}}{({\cal D}^2+z)^{\frac12}}  
\end{array}
$$
define an entire cyclic cocycle \cite{Co}, pp.523,540 over the Banach algebra
$A$ (1.3). 
\end{theorem}

There are similar results for odd Fredholm modules \cite{Co},\,\cite{Co4}.
A few remarks about the formula (1.24) are necessary.
\begin{enumerate}
\item[a)] The expressions
$F(z)\,a^0\,[F(z),\,a^1]\ldots[F(z),\,a^{2n}],\,a^0,\ldots,a^{2n}\in
A,\,z\in\C\backslash\R_-,$ \\
are usually not of trace class or in the domain of
$Trace'$. 
\item[b)]Only the integral described in (1.24) is in the domain of $Trace'$, 
which consists of all bounded operators on $\cal H$ whose even part is of
trace class. 
\item[c)] The formula (1.24) does not explicitely involve the grading
operator $\epsilon$. 
\item[d)] The holomorphic family $F(z),\,z\in\C\backslash\R_-$ is constant
modulo the compact operators: $F(z)\,=\,\frac{{\cal D}}{(1+{\cal
D}^2)^{\frac12}}  \,\,mod\,\,{\cal K(H)}$ and the latter expression does not
involve the grading operator $\epsilon$.
\item[e)] The grading operator $\epsilon$ shows up only in the subdominant 
part $\frac{\epsilon\, z^{\frac12}}{({\cal D}^2+z)^{\frac12}}$ of the
operator family $F(z)$. 
\end{enumerate}
We comment on the grading in such detail because every even Fredholm module 
is homotopic to zero among ungraded modules.

\begin{theorem}\cite{Co}, p.544
Let ${\cal E}\,=\,({\cal H},\,\rho,\,\epsilon,\,{\cal D})$ be an even unbounded 
$\theta$-summable Fredholm module over $A$ and let $e\in K_0(A)$.
Then the index formula
$$
\langle Ch_\epsilon({\cal E}),ch(e)\rangle\,=\,Index({\cal E}_e)
\eqno(1.25)
$$
holds where $ch:K_*(A)\to HC_*^\epsilon(A)$ denotes the Chern character (1.23)
in entire cyclic homology.
\end{theorem}

There is another well known entire cyclic character cocycle
associated to an unbounded $\theta$-summable Fredholm module which
was introduced by Jaffe, Lesniewski and Osterwalder \cite{JLO} and
turns out to be cohomologous to Connes' character cocycle as was shown in
\cite{Co4}. The JLO-cocycle is given by a simpler and more manageable formula 
than the cocycle of Connes. Connes' cocycle can be described in terms of
universal algebras however which will enable us to relate it to the abstract
bivariant  characters. We note finally that there are character formulas for
finitely  summable unbounded Fredholm modules \cite{Co1},\cite{CM1}. They are
cohomologous in entire cyclic cohomology to Connes character cocycle as well
\cite{CM1}.

{\bf The bivariant Chern-Connes character}

\begin{theorem}\cite{Pu1}, 6.3.
There exists a unique natural transformation 
$$
ch_{biv}:\,KK_*(-,-)\,\longrightarrow\,HC_*^{loc}(-,-)
\eqno(1.26)
$$
of bifunctors on the category of separable $C^*$-algebras, which is
multiplicative and not identically zero. It is compatible (up to universal
constants)  with exterior products and excision.
\end{theorem}

The $K$-homology groups of a $C^*$-algebra can be described in terms of
Kasparov theory as $K^*(A)\,=\,KK_*(A,\C)$. Therefore one obtains as special
case of (1.26) a character 
$$
\check{ch}:\,K^*(-)\,=\,KK_*(-,\C)\,\overset{ch_{biv}}{\longrightarrow}\,
HC_*^{loc}(-,\C)\,=\,HC^*_{loc}(-)
\eqno(1.27)
$$
on $K$-homology with values in local cyclic cohomology. 

According to Cuntz \cite{Cu}, \cite{Cu1} the Chern-Connes character on (even)
$K$-homology can be described in terms of universal algebras. We will recall 
a few facts from \cite{Cu}.

Let $A$ be a $C^*$-algebra and let $Q_{C^*}A=A*A$ be the free product of $A$
with itself in the category of $C^*$-algebras. Denote by $\gamma$ the
involutive automorphism of $QA$ switching the two copies of $A$  and let
$\theta:A\to QA$ repectively $\theta^\gamma:A\to QA$ be the canonical
inclusions. There exists a natural extension
$$
0\to q_{C^*}A\to Q_{C^*}A\overset{id*id}{\longrightarrow} A\to 0
\eqno(1.28)
$$
with the two natural multiplicative linear sections $\theta$ and
$\theta^\gamma$. 

Consider the local cyclic (co)homology groups (1.17) associated to this
extension. By excision the natural map 
$$
HC_*^{loc}(A,q_{C^*}A)\,\overset{\simeq}{\longrightarrow}\,Ker\left(\,HC_*^{loc}(A,Q_{C^*}A)\,\to\,HC_*^{loc}(A,A)\right)
\eqno(1.29)
$$
is an isomorphism. We denote by 
$$
\iota^{loc}_A\,\in\,HC_*^{loc}(A,q_{C^*}A)
\eqno(1.30)
$$
the natural element corresponding to the class
$$
(\theta_*-\theta^{\gamma}_*)\,\in\,Ker\left(\,HC_*^{loc}(A,Q_{C^*}A)\,\to\,HC_*^{loc}(A,A)\right)
\eqno(1.31)
$$
under this isomorphism.

Let now ${\cal E}\,=\,({\cal H},\,\rho,\,\epsilon,\,F)$ be an (even) 
bounded Fredholm module over the $C^*$-algebra $A$ which satisfies the
normalization conditions (1.4). There is a canonical pair of homomorphisms
$$
\begin{array}{ccc}
\rho_0,\,\rho_1:A\,\to\,{\cal
L(H)},
& \rho_0(a)\,=\,\frac{1+\epsilon}{2}\,\rho(a)\,\frac{1+\epsilon}{2}, & 
\rho_1(a)\,=\,F\,\frac{1-\epsilon}{2}\,\rho(a)\,\frac{1-\epsilon}{2}\,
F \end{array}  \eqno(1.32)
$$
These homomorphisms coincide modulo compact operators. In particular one
obtains a commutative diagram
$$
\begin{array}{cccclcccc}
 0 & \to & {\cal K(H)} & \longrightarrow & {\cal L(H)} & \longrightarrow &
{\cal Q(H)} & \to & 0 \\
 & & & & & & & & \\
 & & \psi\uparrow & & \uparrow\rho_0*\rho_1 & & \uparrow & & \\
 & & & & & & & & \\
0 & \to & q_{C^*}A & \longrightarrow & Q_{C^*}A & \longrightarrow &
 A & \to & 0 \\
\end{array}
\eqno(1.33)
$$
The map 
$$
\psi:\,q_{C^*}A\,\longrightarrow\,{\cal K(H)}
\eqno(1.34)
$$
is called the characteristic homomorphism associated to $\cal E$.

The local cyclic cohomology of the $C^*$-algebra of compact operators is given
by 
$$
HC^0_{loc}({\cal K(H)})\,\simeq\,\C
\eqno(1.35)
$$
The canonical generator 
$$
[\tau]\,\in\,HC^0_{loc}({\cal K(H)})\,=\,HC_0^{loc}({\cal K(H)},\C)
\eqno(1.36)
$$ 
is cohomologous to the standard trace 
on the ideal $\ell^1({\cal H})$ of trace class operators. 

We exhibit explicit cyclic cocycles 
which represent the generating class $[\tau]$. This is
possible in local cyclic cohomology \cite{Pu}, section 7, but we prefer to
use an analytic cyclic cocycle. We follow therefore \cite{Me},section 3.4.

Let ${\cal P}=(P_n)_{n\in\N}$ be a sequence of finite rank projections 
in $\cal L(H)$ such that 
$\underset{n\to\infty}{lim}\,Rank \,P_n\,=\,\infty,\,\overline{lim}\,(Rank\,
P_n)^{\frac1n}<\infty$. Then the cochain  
$$
\tau_{\cal P}\,=\,Trace\,-\,\partial(\mu_{\cal P}),\,\mu_{\cal
P}=(\mu_{\cal P}^{2n+1})_{n\in\N},
\eqno(1.37)
$$
$$
n!\,\,\mu_{\cal P}^{2n+1}(a^0,\ldots,a^{2n+1})\,=\,
Trace((1-P_n)\, a^0\,(1-P_n)\ldots(1-P_n)\,a^{2n+1}\,(1-P_n))
\eqno(1.38)
$$
defines an analytic cyclic cocycle on the $C^*$-algebra $\cal K(H)$ 
of compact operators on $\cal H$. Its cohomology class does not depend on
the choice of $\cal P$. The image of $\tau_{{\cal P}}$ in local cyclic
cohomology represents the fundamental class $[\tau]\in HC_{loc}^0({\cal
K(H)})$.

With these notations at hand one has the following description of the 
Chern-Connes character on $K$-homology \cite{Cu1}, \cite{Pu}: 

Let ${\cal E}$ be an even bounded Fredholm module over the $C^*$-algebra $A$
and let \\ $[{\cal E}]\in K^0(A)=KK(A,\C)$ be its $K$-homology class.
Then 
$$
\check{ch}([{\cal
E}])\,=\,\iota_A^{loc}\circ\,\psi_*\,\circ\,[tau]
\,\in\,HC_0^{loc}(A,\C)\,=\,HC^0_{loc}(A)
\eqno(1.39) 
$$

In the same spirit it is possible to construct a bivariant character of 
unbounded Kasparov bimodules \cite{BJ}. It depends only on the bivariant
$K$-theory class of the considered bimodule. In \cite{Ni} Nistor introduced a
bivariant Chern-Connes character for finitely summable bimodules over
pre-$C^*$-algebras with values in bivariant periodic cyclic cohomology. It is
compatible with the previously discussed  bivariant character in the
appropriate sense. 

Our main goal will be the comparison of the two characters discussed above. A
first step in this direction is the description of Connes' character 1.5 in the
Cuntz-Quillen picture of cyclic cohomology. This is closely related to, but not
identical with Connes' interpretation of his character in terms of traces on
universal algebras in \cite{Co}.

{\bf Connes' character in the Cuntz-Quillen picture}

Cuntz and Quillen \cite{CQ}, \cite{CQ1} gave a description of the 
periodic cyclic cochain complex of an algebra \cite{Co1} in terms of 
the $X$-complex of the tensor algebra \cite{Qu}. This description
involves some extra structure on the tensor algebra which we recall first.

Let $A$ be an algebra and let $TA$ be the (nonunital) tensor algebra over $A$.
Let $\varrho:A\to TA$ be the canonical linear inclusion and let 
$\omega:A\otimes A\to TA,\, \omega(a,a')=\varrho(a\,
a')-\varrho(a)\,\varrho(a')$ be its curvature \cite{Qu}. 
Let
$$
0\to IA\to TA\overset{\pi}{\longrightarrow} A\to 0
\eqno(1.40)
$$
be the natural extension of algebras with $\varrho$ as natural linear
section. 

There is a natural linear isomorphism \cite{CQ}  
$$
\begin{array}{ccc}
\Omega^{ev} A & \simeq & TA \\
 & & \\
a^0da^1\ldots da^{2n} & \leftrightarrow & \varrho(a^0)\omega(a^1,a^2)
\ldots\omega(a^{2n-1},a^{2n}) \\
 & & \\
\Omega^{odd} A & \simeq & \Omega^1(TA)/[\Omega^1(TA),A] \\
 & & \\
a^0da^1\ldots da^{2n+1} & \leftrightarrow & \varrho(a^0)\omega(a^1,a^2)
\ldots\omega(a^{2n-1},a^{2n})d\varrho(a^{2n+1})\\
\end{array}
\eqno(1.41)
$$
Under this isomorphism the $IA$-adic filtration of $TA$ corresponds to twice 
the degree filtration (Hodge filtration) of $\Omega^{ev}A$. The $IA$-adic
completion of $TA$ will be denoted by $\widehat{TA}$. 

Now we recall Quillen's complex $X^*(A)$. It is the largest
$\Z/2\Z$-graded subcomplex $X^*(A)\subset CC^*_{per}(A)$
of  the periodic cyclic cochain complex consisting of functionals on
$\Omega A$ which vanish on $\Omega^n A$ for $n>1$. Explicitely 
$$
\begin{array}{cc}
X^0(A)\,=\,A', & X^1(A)\,=\,Ker\left(b:\Omega^1A'\to\Omega^2A'\right) \\
\end{array}
\eqno(1.42)
$$
The differentials are given by $b:X^0(A)\to X^1(A),\,d:X^1(A)\to X^0(A)$. Thus
the even and odd cocycles in $X^*(A)$ are given by the traces on $A$
and the closed $A$-bimodule traces on
$\Omega^1A$, respectively.

The results of Cuntz and Quillen which we are going to use are the following

\begin{theorem}\cite{CQ1}
There is a natural diagram of chain homotopy equivalences
$$
\begin{array}{ccccc}
CC^*_{per}(A) & \overset{CC(\pi)}{\longrightarrow} & CC^*_{per}(\widehat{TA}) &
\hookleftarrow & X^*(\widehat{TA}) \\
\end{array}
\eqno(1.43)
$$
\end{theorem}

It is not true however that (1.41) extends to an isomorphism
$CC^*_{per}(A)\sim X^*(\widehat{TA})$ of chain complexes. Instead one has 

\begin{theorem}\cite{CQ1}
Let $A$ be an algebra.
\begin{itemize}
\item[a)] Let $CC^*_{harm}(A)\subset CC^*_{per}(A)$ be the generalized
1-eigenspace of the Karoubi operator (1.11) acting on $CC^*_{per}(A)$. Then 
$CC^*_{harm}(A)$ is a subcomplex and in fact a natural deformation retract 
of $CC^*_{per}(A)$.
\item[b)] Let $X^*(\widehat{TA})_{harm}\subset X^*(\widehat{TA})$ be the image 
of $CC^*_{harm}(A)$ under the linear map (1.41). Then
$X^*(\widehat{TA})_{harm}$ is  a subcomplex and in fact a natural deformation
retract of $X^*(\widehat{TA})$. \item[c)] There is a natural isomorphism of
the harmonic subcomplexes \\ $CC^*_{harm}(A)\simeq X^*(\widehat{TA})_{harm}$
which realizes the  homotopy equivalence (1.43). 
\item[d)] The isomorphism of c) is given on forms of fixed degree by a scalar
multiple of (1.41). 
\end{itemize}
\end{theorem}

In particular (1.41) maps even harmonic cocycles on the cyclic bicomplex
bijectively to harmonic traces on the $IA$-adic completion
$\widehat{TA}$ of the tensor algebra $TA$.

Connes emphasizes in \cite{Co} that his character cocycle $Ch_\epsilon$ is
normalized. This fact implies

\begin{lemma}
The cocycle $Ch_\epsilon$ is invariant under the Karoubi operator $\kappa$
(1.11). In particular it is harmonic in the sense of 1.9. 
\end{lemma}

\begin{proof}
In fact 
$$
\kappa(Ch_{\epsilon}^{2n}(a^0,\ldots,a^{2n}))
=Ch_{\epsilon}^{2n}(a^{2n},a^0,\ldots,a^{2n-1})
-Ch_{\epsilon}^{2n}(1,a^{2n}a^0,a^1,\ldots,a^{2n-1})
$$
$$
=\tau(Fa^{2n}[F,a^0][F,a^1],\ldots,[F,a^{2n-1}])
-\tau(F[F,a^{2n}a^0][F,a^1],\ldots,[F,a^{2n-1}])
$$
$$
=-\tau(F[F,a^{2n}]a^0[F,a^1],\ldots,[F,a^{2n-1}])
=\tau([F,a^{2n}]Fa^0[F,a^1],\ldots,[F,a^{2n-1}])
$$
$$
=\tau(Fa^0[F,a^1],\ldots,[F,a^{2n-1}][F,a^{2n}])=Ch_{\epsilon}^{2n}(a^0,\ldots,a^{2n})
$$
\end{proof} 

The previous lemma allows to identify the cocycle $Ch_\epsilon$
with a trace on the (suitably completed) tensor
algebra $TA$. This makes it possible to give an alternative description of
Connes' character cocycle. It should be noted that our description of the
character cocycle in terms of traces on universal algebras is different from
the one given by Connes in \cite{Co}.

\begin{lemma}
Let $\rho:A\to B$ be a homomorphism of algebras and let $P=P^2\in B$ be an
idempotent element. Let $\varphi:A\to B,\,\varphi(a)=P\rho(a)P$ be the
contraction with $P$ and denote by $T\varphi:TA\to B$ be the corresponding
homomorphism of algebras. Then (in the notations of (1.41))
$$
T\varphi(\varrho(a^0)\omega(a^1,a^2)\ldots\omega(a^{2n-1},a^{2n}))=(-1)^n
P\rho(a^0)P[P,\rho(a^1)]\ldots[P,\rho(a^{2n})] 
\eqno(1.44)
$$
\end{lemma}
\begin{proof}
By definition of the curvature of a linear map (we suppress $\rho$ from the
notation) 
$$
T\varphi(\omega(a,a'))=Paa'P-(PaP)(Pa'P)
$$
$$
=Pa(1-P)a'P=P[a,1-P][a',P]=
-P[P,a][P,a']
$$
Since further $P[P,a][P,a']=P[P,a][P,a']P$
we deduce
$$
T\varphi(\varrho(a^0)\omega(a^1,a^2)\ldots\omega(a^{2n-1},a^{2n}))=(Pa^0P)(-1)^nP[P,a^1]\ldots[P,a^{2n}]
$$
$$
=(-1)^nPa^0P[P,a^1]\ldots[P,a^{2n}]
$$
\end{proof}

\begin{prop}
Let ${\cal E}=({\cal H},\rho,\epsilon,{\cal D})$ be an even unbounded
$\theta$-summable Fredholm module over $A$. Let $F(z)=\frac{{\cal
D}+\epsilon\,z^{\frac12}}{({\cal D}^2+z)^{\frac12}}\in{\cal O}$ be the
corresponding  holomorphic family of operators where $\cal O$ denotes the
algebra of operator valued holomorphic functions on $\C\backslash\R_-$.
Define linear maps
$\varphi_{0,1}:A\to{\cal O}$ by 
$$ 
\begin{array}{lr}
\varphi_0(a)=\frac{1+F(z)}{2}\, \rho(a)\,\frac{1+F(z)}{2}, & 
\varphi_1(a)=\epsilon\, \frac{1-F(z)}{2}\, \rho(a)\,
\frac{1-F(z)}{2}\, \epsilon\\ 
\end{array} 
\eqno(1.45)
$$
and let $T\varphi_{0,1}:TA\to{\cal O}$ be the corresponding homomorphisms of
algebras. Let finally 
$$
\tau(f)\,=\,Trace'\left(\frac{1}{\sqrt{\pi}\,i}\int_{1-i\infty}^{1+i\infty}
f(z^2)\,e^{z^2}\,dz\right),\,\,\,f\in{\cal O}
\eqno(1.46)
$$
be the densely defined unbounded trace on $\cal O$ introduced in (1.24).
Then under the isomorphism of harmonic subcomplexes 1.9 
the Connes character cocycle $Ch_\epsilon({\cal E})\in
CC_{\epsilon}(A)^0_{harm}$ corresponds to the harmonic trace 
$T\varphi_0^*(\tau)-T\varphi_1^*(\tau)\in X^0({\cal T}A)_{harm}$.
$$
\begin{array}{ccc}
CC_{\epsilon}(A)_{harm}^* & \simeq & X^*({\cal T}A)_{harm}\\
 \in & & \ni \\
Ch_\epsilon({\cal E}) & \leftrightarrow &
T\varphi_0^*(\tau)-T\varphi_1^*(\tau) \\
\end{array}
$$
\end{prop}

\begin{proof}
The linear map $\varphi_0$ is given by contraction with the idempotent
$P_+=\frac{1+F(z)}{2}$ and the linear map $\varphi_1$ is the composition of the
contraction with $P_-=\frac{1-F(z)}{2}$ and the conjugation with $\epsilon$.
Lemma 1.11 implies 
$$
T\varphi_0(a^0,\ldots,a^{2n})=(-1)^nP_+a^0P_+[P_+,a^1]\ldots[P_+,a^{2n}]
$$
$$
=(-\frac{1}{4})^n P_+a^0P_+[F,a^1]\ldots[F,a^{2n}]
$$
whereas
$$
T\varphi_1(a^0,\ldots,a^{2n})=
(-\frac{1}{4})^n\epsilon
P_-a^0P_-[F,a^1]\ldots[F,a^{2n}]\epsilon 
$$
Therefore
$$
(T\varphi_0^*-T\varphi_1^*)(\tau)(a^0,\ldots,a^{2n})=
(-\frac{1}{4})^n\tau(P_+a^0P_+[F,a^1]\ldots[F,a^{2n}]-\epsilon
P_-a^0P_-[F,a^1]\ldots[F,a^{2n}]\epsilon)
$$
$$
=(-\frac{1}{4})^n\tau((P_+a^0P_+-P_-a^0P_-)[F,a^1]\ldots[F,a^{2n}])
$$
$$
=(-\frac{1}{4})^n\tau(P_+a^0[F,a^1]\ldots[F,a^{2n}]P_+-P_-a^0[F,a^1]\ldots[F,a^{2n}]P_-)
$$
$$
=(-\frac{1}{4})^n\tau((P_+-P_-)a^0[F,a^1]\ldots[F,a^{2n}])
=(-\frac{1}{4})^n\tau(Fa^0[F,a^1]\ldots[F,a^{2n}])
$$
The isomorphism $CC_{\epsilon}(A)^*_{harm}\simeq X^*({\cal T}A)_{harm}$ 
is given by (a suitable multiple) of (1.41) so that the conclusion follows. 
\end{proof}

\begin{definition}
Let ${\cal E}\,=\,({\cal H}, \rho, \epsilon, {\cal D})$ be an even unbounded 
$\theta$-summable Fredholm module over $A$. Let $\varphi_{0,1}:A\to{\cal O}$
$$ 
\begin{array}{lr}
\varphi_0(a)=\frac{1+F(z)}{2}\, a\, \frac{1+F(z)}{2}, & 
\varphi_1(a)=\epsilon\, \frac{1-F(z)}{2}\, a\, \frac{1-F(z)}{2}\,
\epsilon\\ 
\end{array} 
\eqno(1.47)
$$
be the associated linear maps. Then the algebra homomorphism
$$
\Phi\,=\,T\varphi_0*T\varphi_1:\,Q(TA)=TA*TA\,\longrightarrow\,{\cal O}
\eqno(1.48) 
$$
is called the characteristic homomorphism associated to $\cal E$. 
\end{definition}

\section{Change of regularization in Connes' character formula}

In this section we achieve two things. We show first that the algebraic
description of the class of Connes' character in the Cuntz-Quillen picture,
obtained in the previous section, leads to cocycles which satisfy all required
continuity  and boundedness properties. This is not very helpful yet, because
the use of unbounded operators and the heat kernel regularization allows only
a very limited use of algebric operations, homotopies etc. These are
indispensible however in order to achieve the transgression of Connes'
character to the character of a bounded Fredholm module. The second issue of
this section is therefore a change of regularization in the character formulas.
 
It is well known that the canonical trace $\tau$ on the algebra $\ell^1({\cal
H})$ of trace class operators on Hilbert space does not extend to a trace on 
the $C^*$-closure ${\cal K(H)}$ of all compact operators. The induced map 
$HC^*_{loc}({\cal K(H)})\,\overset{\simeq}{\longrightarrow}\,HC^*_{loc}(\ell^1({\cal
H}))$ on local cyclic cohomology is however an isomorphism \cite{Pu2}.

Therefore the canonical trace, viewed as a local cyclic cocycle, is
cohomologous to a cocycle $\tau_{{\cal P}}$ which extends to all of ${\cal K(H)}$.
Such a cocycle depends on the choice of a bounded approximate unit $\cal P$ in
$\cal K(H)$ and may be analytic \cite{Me}. 

Connes constructs his character cocycle in terms of a superalgebra
$\widetilde{{\cal L}}$ of operator-valued distributions and a canonical odd
trace $\widetilde{\tau}$ on it. We will replace in all formulas his trace 
by a cyclic cocycle $\widetilde{\tau}_{{\cal P}}$ on $\widetilde{{\cal L}}$ which resembles the
regularized cocycles $\tau_{{\cal P}}$ on ${\cal K(H)}$. This will enable us to perform the crucial
transgression in the next section.

Almost the whole section cocnsists of estimates near zero of the Schatten 
norms of operator-valued distributions on the real halfline and of pointwise
estimates of the operator norm of their Laplace transforms. We recommend some
familiarity with the calculations in \cite{Co} in order to digest the
following pages.

{\bf Operator valued distributions and the Laplace transform}

It turns out to be useful to view the holomorphic families 
of operators which occur in Connes' character formula 1.5 as 
the Laplace transforms of certain operator valued distributions. 
This motivated Connes to introduce the following convolution algebra 
of distributions on the positive real halfline. 

\begin{definition} \cite{Co}, p.531.

Let $\cal H$ be a Hilbert space. The convolution algebra 
$\cal L$ of operator valued distributions consists of all tempered
distributions $T$ on the real line satisfying
\begin{itemize}
\item[a)] $Supp(T)\subset[0,\infty[$
\item[b)] There exists a holomorphic operator valued function $t$ 
on some open cone $\bigcup_{s>0}sB(1,r),\,B(1,r)=\{z\in\C,\,\vert
z-1\vert<1\},0<r<1$ in the open right halfplane such that $T$ coincides with
$t$ on $]0,\infty[$
\item[c)] the function 
$$
h(p)=\underset{z\in\frac{1}{p}B(1,r)}{sup}\parallel
t(z)\parallel_p,\,p\in]1,\infty[
$$ 
is of polynomial growth where $\parallel-\parallel_p$ denotes the Schatten
$p$-norm. 
\end{itemize}
\end{definition}

It is a nontrivial fact that $\cal L$ is indeed an algebra (see
\cite{Co},p.533). Connes introduces also the following superalgebra of
distributions.

\begin{definition} \cite{Co}, p.534.

Let $\cal H$ be a Hilbert space and let $\cal L$ be the associated 
algebra 2.1 of distributions. Denote by $\widetilde{{\cal L}}$ the 
quadratic Galois extension of $\cal L$ obtained by adjoining 
a square root $\lambda^{\frac12}$ of the distribution $\lambda=\delta_0'$. 
Then $\lambda^{\frac12}$ is a central element in $\widetilde{{\cal L}}$ and 
$Gal(\widetilde{{\cal L}},{\cal L})\,=\,\Z/2\Z$.
\end{definition}

The Laplace transform is a homomorphism $L$ from the convolution algebra $\cal
L$ of operator valued distributions to the algebra of operator
valued holomorphic functions on the open right half plane
$\{z\in\C,\,Re(z)>0\}$.  It is given by the formula
$$
L(T)(z)\,=\,\int_0^\infty T(s)\, e^{-sz}ds 
\eqno(2.1)
$$
for all $T\in{\cal L}$. We extend it to a homomorphism
$\widetilde{L}$ on $\widetilde{{\cal L}}$ by sending
$\lambda^{\frac12}\in\widetilde{{\cal L}}$ to the branch $z^{\frac12}\in{\cal
O}$ of the square root function which takes the value $+1$ at 1. The Laplace
transform maps $\cal L$ isomorphically onto its picture. An explicit formula 
for the value of the inverse Laplace transform
$\widetilde{L}^{-1}:\,\widetilde{L}(\widetilde{{\cal L}})\to\widetilde{{\cal
L}}$ at $s>0$ is given on the odd part of $\widetilde{{\cal L}}$ by
$$
T(s)\,=\,\frac{1}{2\pi i}
\int_{1-i\infty}^{1+i\infty}S(z)\,
z^{-\frac12}\, e^{s\, z}\, dz\,= \,\int_\Gamma S(z^2)\, d\mu_s(z)
\eqno(2.2)
$$
where $S=\widetilde{L}(\lambda^{\frac12}\, T), T\in{\cal L},$ and 
$$
d\mu_s(z)\,=\,\frac{1}{\pi i}\, e^{sz^2}\, dz 
\eqno(2.3)
$$ 
is a suitable Gaussian measure.
The integration in the last expression is carried out over the oriented curve
$$
\Gamma\,=\,\{z\in\C,\,arg\,(z-1)\,=\,\pm\frac{\pi}{3}\}
\eqno(2.4)
$$
going from the lower to the upper halfplane.
The equality of the two integrals
follows from the Cauchy integral formula  by a reasoning similar to that in
\cite{Co}, pp. 542-543.

We will supress the symbol $\widetilde{L}$ for the Laplace
transform from the notation. If $T\in\widetilde{{\cal L}}$ is an operator
valued distribution, then its Laplace transform will be denoted by
$T(z)$.

\begin{deflem} \cite{Co}, p.534.

Let $\cal H$ be a Hilbert space and let $\widetilde{{\cal L}}$ (see 2.1, 2.2)
be the associated  superalgebra of distributions. Then the functional
$$
\begin{array}{ccc}
\widetilde{\tau}(T)\,=\,Trace(T_-(1)), & T=T_++\lambda^{\frac12}T_-\in
\widetilde{{\cal L}}, & T_\pm\in{\cal L} \\
\end{array}
\eqno(2.5)
$$
defines an (odd) trace on the convolution algebra $\widetilde{{\cal L}}$.
It corresponds under the Laplace transform to a multiple of the trace
functional introduced in (1.24).
\end{deflem}

\begin{definition}
Let ${\cal P}=(P_n)_{n\in\N}$ be a sequence of finite rank projections 
in $\cal L(H)$ such that 
$\underset{n\to\infty}{lim}\,Rank \,P_n\,=\,\infty,\,\overline{lim}\,(Rank\,
P_n)^{\frac1n}<\infty$. Then the functional cochain  
$$
\widetilde{\tau}_{\cal
P}\,=\,\widetilde{\tau}\,-\,\partial(\widetilde{\mu}_{\cal
P}),\,\widetilde{\mu}_{\cal P}=(\mu_{\cal P}^{2n+1})_{n\in\N}, \eqno(2.6) 
$$
$$
n!\,\,\widetilde{\mu}_{\cal P}^{2n+1}(T^0,\ldots,T^{2n+1})\,=\,
\widetilde{\tau}((1-P_n)\, T^0\,(1-P_n)\ldots(1-P_n)\,T^{2n+1}\,(1-P_n))
$$
is a cocycle on the chain complex $CC^0(\widetilde{{\cal L}})$ which is
cohomologous to the canonical trace $\widetilde{\tau}$ on $\widetilde{{\cal
L}}$ 
\end{definition}

Connes associates to every unbounded $\theta$-summable Fredholm
module a couple of operator valued distributions. These are used to establish
the continuity of the character cocycle (1.24).

\begin{deflem} \cite{Co},p.535.
Let ${\cal E}\,=\,({\cal H}, \rho, \epsilon, {\cal D} )$ 
be an even unbounded $\theta$-summable Fredholm module over $A$.
\begin{itemize}
\item[i)]
There exist operator-valued distributions \cite{Co},p.535,
$$
\begin{array}{ccc}
F\,=\,{\cal D}N\,+\,\lambda^{\frac12}\,\epsilon\, N\,\in\widetilde{{\cal
L}}, & N\,=\,\frac{1}{\sqrt{\pi s}}\, e^{-s{\cal D}^2}\,\in\,{\cal L}, &
\lambda^{\frac12}*\lambda^{\frac12}=\delta_0' \\ 
\end{array} \eqno(2.7)
$$
satisfying 
$$
F^2\,=\,\delta_0\,=\,1\,\eqno(2.8)
$$
\item[ii)]
The Laplace transform of the distribution in i) is given by 
$$
\widetilde{L}(F)(z)\,=\,\frac{{\cal D}\,+\,z^{\frac12}\,\epsilon}{({\cal
D}^2\,+\,z)^{\frac12}}
\eqno(2.9)
$$
\end{itemize}
\end{deflem}

\begin{definition}
We denote by $\varphi_0,\,\varphi_1:A\to\widetilde{{\cal L}}$ the linear maps
$$
\begin{array}{cc}
\varphi_0(a)\,=\,\frac{1+F}{2}\,\rho(a)\,\frac{1+F}{2}, & 
\varphi_1(a)\,=\,\epsilon\,\frac{1-F}{2}\,\rho(a)\,\frac{1-F}{2}\,\epsilon
\\ \end{array}
\eqno(2.10)
$$
whose Laplace transform yields the maps introduced in 1.13. 
We call the induced algebra homomorphism
$$
\Phi=T\varphi_0*T\varphi_1:Q(TA)\to\widetilde{{\cal L}}
\eqno(2.11)
$$ 
still the characteristic morphism associated to ${\cal E}$. Its Laplace
transform is the morphism of 1.13. 
\end{definition}

{\bf Majorizing functions}

We recall a few facts from \cite{Co}.

\begin{definition}\cite{Co}, page 537.

Let $B(1,r)\,=\,\{z\in\C,\vert z-1\vert<r\}$. Then $f\in L^1([0,1])$ is a
majorizing function for $T\in\cal L$ with respect to $B(1,r)$  if
\begin{itemize}
\item[a)] $\underset{z\in\frac{1}{p} B(1,r)}{\sup}\parallel
T(z)\parallel_p\,\leq\,f(\frac{1}{p}),\,\forall p\in[1,\infty]$ \\
\item[b)] $T(\Phi)\,=\,\int_0^\infty T(s)\Phi(s) ds,\,\forall\Phi\in{\cal
S}(\R)$ \\ 
\end{itemize} 
We write $T \underset{r}{\prec} f $ in this case.
\end{definition}

Remark: Not every distribution in $\cal L$ possesses a majorizing function.

\begin{lemma}\cite{Co}

Let $f\in L^1([0,1])$ be a
majorizing function for $T\in\cal L$ with respect to $B(1,r)$ such that 
$s^{-k}\, f\in L^1([0,1])$. Then
$$
\begin{array}{ccc}
\frac{\partial^k T}{\partial s^k} & \underset{\frac{r}{2}}{\prec} &
\left(\frac{2}{r}\right)^k\, k!\, (s^{-k}f) \\
\end{array} \eqno(2.12)
$$
\end{lemma}

\begin{lemma}\cite{Co}, Lemma 1, page 537.
$$
\begin{array}{ccccc}
T_0 \underset{r}{\prec} f_0, & T_1 \underset{r}{\prec} f_1 & \Longrightarrow &
T_0*T_1 \underset{r}{\prec} 2f_0*f_1, & \forall r<1\\
\end{array} \eqno(2.13)
$$
\end{lemma}

\begin{lemma}
Define smooth functions $f_{n}$ on $\R_+$ for $n\geq 1$ by 
$$
\begin{array}{cc}
 f_{1}(s)=\frac{1}{\sqrt{\pi s}}, & 
f_{n+m}=f_{n}*f_{m} \\
\end{array} 
$$ 
where the product is given by convolution 
on the additive semigroup $(\R_+,+)$. Then
$$
\begin{array}{cc}
f_{2n}(s)=\frac{s^{n-1}}{(n-1)!}, &
f_{2n+1}(s)=\frac{1}{\sqrt{\pi}}\, 4^{n}\,\frac{(n)!}{(2n)!}\,
s^{n-\frac{1}{2}} \\
\end{array}
\eqno(2.14)
$$
\end{lemma}
\begin{proof}
Take the Laplace transforms $\widehat{f}(s)=\int_0^\infty f(t)e^{-st}dt$ 
of the functions in question.
\end{proof}

\begin{lemma}
Let $T\in{\cal L}$ be an operator valued distribution, let $0<r<1$ and let 
$n>0$ be an integer. Then
$$
T\,\underset{r}{\prec}\,f_{n+2k}(s)
\,\Rightarrow\,
\frac{\partial^k T}{\partial s^k}\,\underset{\frac{r}{2}}{\prec}\,C(r)^k\,
f_{n}(s) \eqno(2.15)
$$
\end{lemma}
This is clear from 2.8, 2.9 and 2.10.

{\bf Estimates of Schatten norms of operator valued distributions}

We fix now an even unbounded $\theta$-summable Fredholm module 
${\cal E}\,=\,({\cal H}, \rho, \epsilon, {\cal D})$
 over $A$ and recall that $A$ is supposed to be a Banach algebra with respect
to the norm
$$
\parallel a\parallel_A\,=\,\parallel \rho(a)\parallel_{{\cal
L(H)}}\,+\, \parallel[{\cal D}, \rho(a)]\parallel_{{\cal
L(H)}},\,\,\,a\in{\cal A}
$$ 
We will omit $\rho$ from the notation if no
confusion is likely to arise.

Following Connes we consider the
associated operator-valued distributions (2.7).

\begin{lemma}\cite{Co}, pp.535, 537.
Let ${\cal E}=({\cal H},\rho, {\cal D},\epsilon)$ be an even unbounded $\theta$-summable
unbounded Fredholm module over $A$ and let $N,\,{\cal D}N$ be the associated
distributions. Let $P$ be any spectral projection associated to ${\cal D}^2$
and put  $\eta\,=\,Trace(P\, e^{-\frac{{\cal D}^2}{4}})$.
\begin{itemize}
\item[a)] 
$$
\begin{array}{cc}
P\,{\cal D}N \,=\,PS\,\delta_0\,+\,(1-\lambda)\, PT, & 
S\in{\cal L(H)},T\in{\cal L}\\
\end{array} \eqno(2.16)
$$  
with 
$$
\begin{array}{cc}
PT\underset{\frac14}{\prec}C_1({\cal
E})\,\eta^s\, f_{1}(s), & 
PN\underset{\frac14}{\prec}C_2({\cal E})\,\eta^s\, f_{1}(s) 
\end{array} 
\eqno(2.17) 
$$ 
(\cite{Co}, Lemma 2 b) p. 535.

\item[b)]
For all $a\in A$ one
has 
$$
\begin{array}{cc}
P[F,a]P\,=\,P[F,a]_+
P\,+\,\lambda^{\frac12}\,P[F,a]_-P, & [F,a]_\pm\in{\cal L} \\
\end{array} 
$$
with
$$
\begin{array}{cc}
P[F,a]_+P\,\underset{\frac14}{\prec}\,C_3({\cal
E})\eta^s f_1(s)\parallel a\parallel_A, &  
P[F,a]_-P\,\underset{\frac14}{\prec}\,C_4({\cal
E})\eta^s f_2(s)\parallel a\parallel_A  \\ 
\end{array}
\eqno(2.18)  
$$
\item[c)] Suppose that $Im(P)\cap Ker({\cal D})=0$ and that $\eta < 1$.
Then  
$$
PT\underset{\frac14}{\prec}C_5({\cal
E})(-4\,log\,\eta)^{-1}\eta^s f_{1}(s). 
\eqno(2.19) 
$$ 
\end{itemize}
\end{lemma}

\begin{proof}
This is a mild generalization of \cite{Co}, Lemma 2, p.537. The proof of a)
and b) is exactly the same as in loc. cit. Concerning c) we observe (in the
notations of \cite{Co}) that
$$
PT\,=\,\frac{1}{\pi}\int_0^\infty
P{\cal D}({\cal D}^2+1+\rho)^{-1}(P e^{-s({\cal
D}^2+\rho)})\rho^{-\frac12}d\rho 
$$
Let $\lambda_1>0$ be the smallest eigenvalue of the restriction of ${\cal
D}^2$ to the image of $P$. Then 
$$
\parallel P{\cal D}({\cal
D}^2+1+\rho)^{-1}\parallel\,\leq\,\lambda_1^{-1}\, \leq\,(-4\,log\,\eta)^{-1}
$$
because $e^{-\frac{\lambda_1}{4}}\,<\,Trace(P\,e^{-\frac{{\cal
D}^2}{4}})\,=\,\eta$. On the other hand one finds for the Schatten norm
$$
\parallel\frac{1}{\pi}\int_0^\infty (P\,e^{-s({\cal
D}^2+\rho)})\rho^{-\frac12}d\rho\parallel_{\frac{1}{s}}\,\leq\,
\eta^s\frac{1}{\pi}\int_0^\infty e^{-s\rho}\rho^{-\frac12}d\rho
\,=\,\eta^s f_1(s) 
$$
which suffices to establish our claim.
\end{proof}

Let $TA$ be the (nonunital) tensor algebra over $A$. Let $\varrho:A\to TA$
be the canonical linear inclusion and let $\omega:A\otimes A\to TA,\,
\omega(a,a')=\varrho(a\, a')-\varrho(a)\,\varrho(a')$ be its
curvature \cite{Qu}. Recall that there is a natural linear isomorphism
$$
\begin{array}{ccc}
\Omega^{ev} A & \simeq & TA \\
 & & \\
a^0da^1\ldots da^{2n} & \leftrightarrow & \varrho(a^0)\omega(a^1,a^2)
\ldots\omega(a^{2n-1},a^{2n}) \\
\end{array}
$$

\begin{remark}
 In the sequel we only will prove estimates for tensors of the 
form $\varrho\omega^n$ and not of the form $\omega^n$, which correspond to
exact forms under (1.41). This is purely a matter of notational convenience.
All estimates of this section are valid for tensors of the latter
form as well. For them the proofs are usually even simpler.
\end{remark}

\begin{lemma}
Let $S\subset A$ be a bounded set and put $\parallel
S\parallel\,=\,\underset{a\in S}{Sup}\parallel a\parallel$. 
Let $\alpha\in TA$ be an element of the form
$\alpha\,=\,\varrho\omega^n(a^0,\ldots,a^{2n}),\,a^0,\ldots a^{2n}\in S,\,n>0$.
Let $T\varphi_0,
T\varphi_1:TA\to\widetilde{{\cal L}}$ be the homomorphisms introduced  in 2.6. 
Then 
$$
\begin{array}{cc}
T\varphi_{i}(\alpha)\,=\,T\varphi_i(\alpha)_+\,+\,\lambda^{\frac12}\,
T\varphi_{i}(\alpha)_-, &  T\varphi_{i}(\alpha)_\pm\in{\cal L},i=0,1, \\
\end{array}
$$
with
$$
\begin{array}{cc}
T\varphi_{i}(\alpha)_+\,\underset{\frac{1}{16}}{\prec}\,(C_6({\cal
E})\parallel S\parallel)^{2n+1}
f_{2n-1}(s), &
T\varphi_{i}(\alpha)_-\,\underset{\frac{1}{16}}{\prec}\,
(C_7({\cal E})\parallel S\parallel)^{2n+1}
\, f_{2n}(s)\\ 
\end{array} 
\eqno(2.20)
$$
\end{lemma}

\begin{proof}
This is essentially Lemma 3 on page 538 of \cite{Co}. As similar arguments
will be used frequently in this section we will give all necessary details. 

The calculations in 1.12 show that
$$
T\varphi_{0}(\varrho(a^0)\omega(a^1,a^2)\ldots\omega(a^{2n-1},a^{2n}))\,=\,
(-\frac{1}{4})^n P_+a^0P_+[F,a^1]\ldots[F,a^{2n}]
$$
$$
T\varphi_{1}(\varrho(a^0)\omega(a^1,a^2)\ldots\omega(a^{2n-1},a^{2n}))\,=\,
(-\frac{1}{4})^n \epsilon P_-a^0P_-[F,a^1]\ldots[F,a^{2n}]\epsilon
$$
Now
$$
[F,a^1]\ldots[F,a^{2n}]=T_+\,+\,\lambda^{\frac{1}{2}}T_-,\,\,T_+,T_-\in{\cal
L},  
$$ 
with
$$
\begin{array}{cc}
T_+=\sum_{j=0}^n\underset{J'\subset\{1,\ldots,2n\},\vert
J'\vert=2j}{\sum}\lambda^j\,T_{J'}, &
T_-=\sum_{j=0}^{n-1}\underset{J''\subset\{1,\ldots,2n\},\vert
J''\vert=2j+1}{\sum}\lambda^j\,T_{J''} \\
\end{array}
$$
where for $J\subset\{1,\ldots,2n\}$ one has
$T_J\,=\,S_1\ldots S_{2n}$ with $S_i=[F,a^i]_-$ if $i\in J$ and
$S_i=[F,a^i]_+$ otherwise. From this one derives with the help of 2.8, 2.9,
2.11 and 2.12 the estimates
$$
T_{J'}\,\underset{\frac14}{\prec}\,2^{2n-1} C_8^{2n-2j}
C_9^{2j}\parallel a^1\parallel_A\ldots\parallel
a^{2n}\parallel_A f_{2n+2j}(s) 
$$
$$
\leq(2C_8+2C_9)^{2n}\parallel
S\parallel^{2n} f_{2n+2j}(s)
$$
and
$$
\lambda^j\,T_{J'}\,\underset{\frac18}{\prec}\,
C_{10}^{2n}\parallel
S\parallel^{2n} f_{2n}(s)
$$
Thus
$$
T_+(a^1,\ldots,a^{2n})\,\underset{\frac18}{\prec}\,(C_{11}({\cal E})\parallel
S\parallel)^{2n} f_{2n}(s)
$$
A similar calculation shows
$$
T_-(a^1,\ldots,a^{2n})\,\underset{\frac18}{\prec}\,(C_{12}({\cal E})\parallel
S\parallel)^{2n} f_{2n+1}(s)
$$
We consider now the distributions $P_\pm a^0P_\pm$. One finds from 2.5,
2.6, and 2.12 
$$
P_\pm a
P_\pm\,=\,P_\pm a
\,-\,P_\pm[P_\pm,
a]\,=\,
S^\pm(a)\,\delta_0\,+\,T_1^\pm(a)\,+\,\lambda^{\frac12}\,
T_2^\pm(a)\,+\,\lambda\, T_3^\pm(a) 
\eqno(2.21)
$$ 
$$
\begin{array}{ccc}
\parallel S(P_\pm)(a)\parallel_{{\cal L(H)}}\leq
C_{13}\parallel a\parallel_A, &
T_i^\pm(a)\,\underset{\frac18}{\prec}\, C_{14}\parallel
a\parallel_A f_{1}(s), & i=1,2,3\\
\end{array}
\eqno(2.22)
$$
Altogether 
$$
T\varphi_{0}(\varrho(a^0)\omega(a^1,a^2)\ldots\omega(a^{2n-1},a^{2n}))\,=
$$
$$
=(S^+\delta_0\,+\,T_1^+\,+\,\lambda^{\frac12}\,
T_2^+\,+\,\lambda\,T_3^+)(a^0)(T_+(a^1,\ldots,a^{2n})\,+\,\lambda^{\frac12}\,
T_-(a^1,\ldots,a^{2n}))
$$
$$
=T\varphi_{0}(\varrho\omega^n)_+\,+
\,\lambda^{\frac12}\,
T\varphi_{0}(\varrho\omega^n)_-
$$
where according to 2.9, 2.10, 2.11,
$$
T\varphi_{0}(\varrho\omega^n)_+ 
\,\underset{\frac{1}{16}}{\prec}\,( C_{15}({\cal E})\parallel
S\parallel)^{2n+1} f_{2n-1}(s)
$$
$$
T\varphi_{0}(\varrho\omega^n)_- 
\,\underset{\frac{1}{16}}{\prec}\,( C_{16}({\cal E})\parallel
S\parallel)^{2n+1} f_{2n}(s)
$$
A similar calculation applies for $T\varphi_{1}$.
\end{proof}

Let $A$ be an algebra and let $QA=A*A$ be the free product of $A$ with itself.
Denote by $\gamma$ the involutive automorphism switching the two copies of $A$ 
and let $\theta:A\to QA$ be the inclusion of the first copy of $A$ into the
free product. For $a\in A$ denote by $pa$ and $qa$ the even and
the odd parts of $\theta(a)$ with respect to the involution $\gamma$. 
Recall that there is a natural linear isomorphism \cite{CQ}, 1.3
$$
\begin{array}{ccc}
\Omega A & \overset{\simeq}{\longrightarrow} & QA \\
 & & \\
a^0da^1\ldots da^n & \leftrightarrow & pa^0qa^1\ldots qa^n \\
\end{array}
\eqno(2.23)
$$

\begin{remark}
For simplicity we will only prove estimates for elements in a free
product of the form $pq^n$ and not of the form $q^n$
corresponding to exact forms under (2.23). All estimates of this section 
are however valid for elements of the latter form and 
for them the proofs are even simpler. 
\end{remark}

\begin{lemma}
Let $\Phi:Q(TA)\to\widetilde{{\cal
L}}$ be the characteristic homomorphism (2.10). Let $S\subset A$ be a bounded
subset of $A$. 
\begin{itemize}
\item [a)] For every element $\alpha\in Q(TA)$ of the form 
$$
\alpha\,=\,p(\varrho\omega^{n_0})q(\varrho\omega^{n_1})\ldots
q(\varrho\omega^{n_k})
(a^0,\ldots,a^{k+2l}),
\eqno(2.24)
$$
$$
a^0,\ldots,a^{k+2l}\in
S,\,l=(n_0+\ldots+n_k) > 0,\,k\geq 0,\,m=\sharp\{i,\,n_i>0\}\leq l  
$$
one has 
$$
\Phi(\alpha)\,=\,\sum_{i=0}^2\lambda^{\frac{i}{2}}\,\Phi_i(\alpha),\,\Phi_i(\alpha)\in{\cal L},
\eqno(2.25)
$$
with
$$ 
\begin{array}{cccc}
\Phi_i(\alpha) & \underset{\frac{1}{64}}{\prec} & ( C_{17}({\cal
E})\parallel S\parallel)^{k+2l+1}
f_{2l-m}(s), & i=0,1 \\  
\Phi_2(\alpha) & \underset{\frac{1}{64}}{\prec} & ( C_{18}({\cal
E})\parallel S\parallel)^{k+2l+1}
f_{2l-m+1}(s) & \\
\end{array}\eqno(2.26)
$$
\item[b)] For every element 
$$
\alpha\,=\,p(\varrho)q(\varrho)^k(a^0,\ldots,a^{k}),
\eqno(2.27)
$$
$a^0,\ldots,a^{k}\in S,\,k\geq 0,$ one has 
$$
\Phi(\alpha)\,=\,S(\alpha)\,\delta_0\,+\,
\sum_{j=0}^3\lambda^{\frac{j}{2}}\,\Phi_j(\alpha),\,\Phi_j(\alpha)\in{\cal
L},
\eqno(2.28)
$$
with
$$
\begin{array}{cccc}
\Phi_j(\alpha) & \underset{\frac{1}{64}}{\prec} & ( C_{19}({\cal
E})\parallel S\parallel)^{k+1}
f_{1}(s), & j=0,1,2 \\  
\Phi_3(\alpha) & \underset{\frac{1}{64}}{\prec} & ( C_{20}({\cal
E})\parallel S\parallel)^{k+1}
f_{2}(s) & \\
\end{array}\eqno(2.29)
$$
The term $S(\alpha)$ occurs only if $k=0$ and then $\parallel
S(\alpha)\parallel\,\leq\,C_{21}({\cal E})\parallel S\parallel$.
\end{itemize}
\end{lemma}

\begin{proof}

If $\alpha=q(\varrho\omega^n(a^0,\ldots,a^{2n})),\,a^0,\ldots,a^{2n}\in
S,\,n>0,$ then  
$$
\Phi(\alpha)_+\,=\,\frac12(T\varphi_0(\alpha)_+-T\varphi_1(\alpha)_+)
 \underset{\frac{1}{16}}{\prec} (C_{22}({\cal
E})\parallel S\parallel)^{2n+1} f_{2n-1}(s) 
$$
$$
\Phi(\alpha)_-\,=\,\frac12(T\varphi_0(\alpha)_--T\varphi_1(\alpha)_-)
 \underset{\frac{1}{16}}{\prec} (C_{23}({\cal
E})\parallel S\parallel)^{2n+1} f_{2n}(s)
$$
by lemma 2.14. A similar estimate holds for
$\alpha=p(\varrho\omega^n(a^0,\ldots,a^{2n}))$,
$a^0,\ldots,a^{2n}\in S,\,n>0.$
If $\alpha=q(\varrho(a)),\,a\in S,$ we find the following.
$$
\Phi(q\varrho(a))\,=\,\frac12(P_+ a P_+-\epsilon
P_- a P_-\epsilon)
$$
$$
=\frac12(P_+\epsilon  a(\epsilon P_+-P_-\epsilon)\,+\,(P_+\epsilon-\epsilon
P_-) a P_-\epsilon)=
$$
$$
\frac12\left(\frac{1+F}{2}\epsilon a\lambda^{\frac12}N\,+\,
\lambda^{\frac12}N a\frac{1-F}{2}\epsilon\right)
$$
$$
=\frac14\lambda^{\frac12}\epsilon\left((N a+a N)\,-\,N[{\cal D},a]
N\right) \eqno(2.30)
$$
Thus 
$$
\Phi(q\varrho(a))\,=\,\lambda^{\frac12}\, T_0(a),
\,\,\,T_0(a)\underset{\frac{1}{4}}{\prec} C_{24}\parallel
a\parallel_A f_{1}(s) 
\eqno(2.31)
$$
by 2.9 and 2.12. A calculation using 2.9 and 2.12 again shows that
$$
\Phi(p(\varrho(a)))\,=\,S(a)\,\delta_0+\sum_{j=0}^3\lambda^{\frac{j}{2}}\,\Phi_j,
\eqno(2.32)
$$
$$
\parallel S\parallel\leq  C_{25}\parallel S\parallel,\,
\Phi_j\underset{\frac{1}{32}}{\prec}  C_{25}\parallel S\parallel
f_1(s),0\leq j\leq 2,\,\Phi_3\underset{\frac{1}{32}}{\prec}  C_{26}\parallel
S\parallel f_3(s)
\eqno(2.33)
$$
Products of the form
$p(\varrho\omega^{n_0})q(\varrho\omega^{n_1})\ldots
q(\varrho\omega^{n_k})$ of the terms considered up to now can be analyzed in
the same way as in the proof of lemma 2.14.  
\end{proof}

\begin{lemma}
Let $P$ be a spectral projection associated to ${\cal
D}^2$ and put $\eta\,=\,Trace(P\, e^{-\frac{{\cal D}^2}{4}})$. 
Define linear maps $(T\varphi_i)^P:TA\to{\widetilde{\cal L}}$ by  
$$
(T\varphi_0)^P(\varrho\omega^n(a^0,\ldots,a^{2n}))\,=\,(-\frac14)^n
P(P_+a^0P_+)P[F,a^1]P\ldots P[F,a^{2n}]P
$$
$$
(T\varphi_1)^P(\varrho\omega^n(a^0,\ldots,a^{2n}))\,=\,(-\frac14)^n
\epsilon P(P_-a^0P_-)P[F,a^1]P\ldots P[F,a^{2n}]P\epsilon
$$
and put
$p^P\,=\,\frac12((T\varphi_0)^P+(T\varphi_1)^P),\,q^P\,=\,\frac12((T\varphi_0)^P-(T\varphi_1)^P)$.
Let finally \\ $\Phi^P:Q(TA)\to\widetilde{{\cal L}}$ be the linear map
defined by  
$$
\Phi^P(p(\varrho\omega^{n_0})q(\varrho\omega^{n_1})\ldots
q(\varrho\omega^{n_k}))\,=\,p^P(\varrho\omega^{n_0})
q^P(\varrho\omega^{n_1})\ldots q^P(\varrho\omega^{n_k})
$$
Then the statement of lemma 2.16 holds for $\Phi^P$ (instead of $\Phi$) 
with the following two modifications: the constants
$C_{17}({\cal E}),\ldots,C_{21}({\cal E})$  have possibly to be changed and all
majorizing functions $\Phi_i,\,\Phi_j$  have to be replaced by
$\eta^s\Phi_i,\,\eta^s\Phi_j$.

Suppose moreover that $Im(P)\cap ker({\cal D})=0$ and that
$\eta\,<\,1$. 
Let $S\subset A$ be a bounded subset of $A$. Then

\begin{itemize}
\item[a)]
For all elements $\alpha\in q(TA)$ of the form (2.24) 
one has 
$$
\Phi^P(\alpha)\,=\,\sum_{i=0}^3\lambda^{\frac{i}{2}}\,\Phi^P_i(\alpha),\,\Phi^P_i(\alpha)\in{\cal L},
\eqno(2.34)
$$
with
$$ 
\begin{array}{cccc}
\Phi^P_i(\alpha) & \underset{\frac{1}{64}}{\prec} & ( C_{27}({\cal
E})\parallel S\parallel)^{k+2l+1}
f_{2l-m+i}(s), & i=0,1 \\  
\Phi^P_i(\alpha) & \underset{\frac{1}{64}}{\prec} & ( C_{28}({\cal
E})\parallel S\parallel)^{k+2l+1} (-4\,log\,\eta)^{-1}
f_{2l-m-1+i}(s) & i=2,3\\
\end{array}\eqno(2.35)
$$
\item[b)] 
For all elements $\alpha\in q(TA)$ of the form (2.27) 
one has 
$$
\Phi^P(\alpha)\,=\,
\sum_{j=0}^3\lambda^{\frac{j}{2}}\,\Phi^P_j(\alpha),\,\Phi^P_j(\alpha)\in{\cal
L},
\eqno(2.36)
$$
with 
$$ 
\begin{array}{cccc}
\Phi^P_j(\alpha) & \underset{\frac{1}{64}}{\prec} & ( C_{29}({\cal
E})\parallel S\parallel)^{k+1}
f_{1}(s), & j=0,1 \\  
\Phi^P_2(\alpha) & \underset{\frac{1}{64}}{\prec} & ( C_{30}({\cal
E})\parallel S\parallel)^{k+1}
((-4\,log\,\eta)^{-1} f_1(s)\,+\,f_{2}(s)), & \\
\Phi^P_3(\alpha) & \underset{\frac{1}{64}}{\prec} & ( C_{31}({\cal
E})\parallel S\parallel)^{k+1}
(-4\,log\,\eta)^{-1} f_{2}(s) & \\
\end{array}\eqno(2.37)
$$
\end{itemize}
\end{lemma}

The proof is similar to that of 2.16. The only difference is that in a) and b)
the estimate (2.17) may be replaced by the sharper bound (2.19).

\begin{lemma}
The notation of 2.15 and 2.16 is understood.
Let $\alpha_0,\ldots,\alpha_{2n+1}\in q(TA)$ be elements of the form 
$$
\alpha_i\,=\,p(\varrho\omega^{n_0})q(\varrho\omega^{n_1})\ldots
q(\varrho\omega^{n_{k_i}})(a^0,\ldots,a^{k_i+2l_i}),
$$
$$
a^0,\ldots,a^{k_i+2l_i}\in
S,\,l_i=(n_0+\ldots+n_{k_i})\geq 0,\,k_i>0,
$$
and put
$k=k_0+\ldots+k_{2n+1},\,l=l_0+\ldots+l_{2n+1}$. Let
$P$ be a spectral projection associated to the self adjoint operator  ${\cal
D}^2$ and put $\eta\,=\,Trace(P\, e^{-\frac14{\cal D}^2})$. Then the odd
part of the operator valued distribution   
$$
T^{\mu_{2n+1}}\,=\,T^{\mu_{2n+1}}_+\,+\,\lambda^{\frac12}\,T^{\mu_{2n+1}}_-\,=
\,P\Phi^P(\alpha_0)P\ldots
P\Phi^P(\alpha_{2n+1}) P 
\eqno(2.38)
$$
can be written as a sum 
$$
T^{\mu_{2n+1}}_-\,=\,\sum_{j=0}^{n+1}\lambda^j\, T_{j}^{\mu_{2n+1}}
\eqno(2.39)
$$
of distributions which satisfy 
$$
T^{\mu_{2n+1}}_{j}  \underset{\frac{1}{128}}{\prec} (C_{32}({\cal
E})\parallel S\parallel)^{k+2l+2n+2}\,f_{l+1}(s)
\eqno(2.40)
$$
Suppose in addition that $Im(P)\cap Ker({\cal D})=0$ and that
$\eta<e^{-\frac14}$.  Then
$$
T^{\mu_{2n+1}}_{j}  \underset{\frac{1}{128}}{\prec} (C_{33}({\cal
E})\parallel
S\parallel)^{k+2l+2n+2}\,(-4\,log\,\eta)^{-j}\,\eta^s\,
f_{l+1}(s)  
\eqno(2.41)
$$
\end{lemma} 

This follows in a tedious but straightforward calculation from 2.9, 2.16 and
2.17. It should be noted that the estimates in (2.40) are sufficient for our
purpose but far from being sharp. Lemma 2.16 and 2.17 are much more precise but
we prefer to give a very simple final formula in (2.40).

{\bf Free products}

Before we proceed we recall a few facts about locally convex
topologies on the various universal algebras we are going to use.

It is easy to see that colimits (like free products for example) 
do not necessarily exist in the category of Banach algebras. This obstacle 
disappears if one passes to a smaller category,
the category of Banach algebras and
contractive homomorphisms (i.e. homomorphisms of norm less or equal too one) as
morphisms.  We do not discuss this in full generality but treat only the
cases we are interested in. 

\begin{definition}
Let $A$ be a Banach algebra let $R\geq 1$ be a real
number. We denote by $Q_RA$ the free product (in the categry of Banach
algebras and contractive homomorphisms) of two copies of $A$ equipped with the
rescaled norm $R\parallel-\parallel_A$. 
\end{definition}

The Banach algebra $Q_RA$ is thus the completion of the algebraic free product
$QA=A*A$ with respect to the largest submultiplicative seminorm for which 
the canonical inclusions $\theta:A\to QA,\,\theta^\gamma:A\to QA$ satisfy
$$
\begin{array}{cc}
\parallel\theta\parallel\leq R, &
\parallel\theta^\gamma\parallel\leq R
\end{array}
$$
By construction $Q_RA$ possesses the following universal property:
there is a canonical
bijection between the set  of pairs of homomorphisms of norm less or equal to
$R$ from $A$ to some Banach algebra $B$ and the set of contractive
homomorphisms (i.e. homomorphisms of  norm less or equal to 1) from $Q_RA$ to
$B$. There are canonical  contractive homomorphisms $Q_{R'}A\to Q_RA$ for
$R'>R$ and the inverse limit 
$$
Q_{top}A\,=\,\underset{\infty\leftarrow R}{\lim} Q_RA
\eqno(2.42)
$$
is a Fr\'echet algebra. It contains the algebraic free product $QA$ as dense
subalgebra and possesses an obvious universal property.

There exists an extension 
$$
0\,\longrightarrow\,q_RA\,\longrightarrow\,Q_RA\overset{id*id}{\longrightarrow}\,A\,\longrightarrow\,0
\eqno(2.43)
$$
with canonical linear multiplicative sections $\theta$ and $\theta^\gamma$ of
norm less or equal to $R$. It is natural in an appropriate sense. Similarly
there is a splitting extension 
$$
0\,\longrightarrow\,q_{top}A\,\longrightarrow\,Q_{top}A\overset{id*id}{\longrightarrow}\,A\,\longrightarrow\,0
$$
of Fr\'echet algebras.
We denote by 
$$
\begin{array}{ccc}
\iota_A^\alpha\,\in\,HC_0^\alpha(A,\,q_RA), & R\geq 1, &
\alpha\,=\,\epsilon,\,anal,\,loc \\ 
\end{array}
\eqno(2.44)
$$
the natural element corresponding to the class
$$
\begin{array}{ccc}
\theta_*-\theta^\gamma_* & \in & Ker(HC_0^\alpha(A,\,Q_RA)\to
HC_0^\alpha(A,\,A)) 
\end{array}
\eqno(2.45)
$$
under the excision isomorphism
in entire, analytic, or local cyclic cohomology, respectively.

We pass now to locally convex topologies on tensor algebras. 
\begin{definition}\cite{Pu2}.
Let $A$ be a Banach algebra and let $R\geq 1$ be a real
number. The Banach algebra $T_{R}A$ is the completion 
of the tensor algebra $TA$ over $A$ with respect to the 
largest submultiplicative seminorm for which the canonical 
linear inclusion $\varrho:A\to TA$ and the bilinear 
curvature map $\omega:A\otimes_\pi A\to TA$ satisfy
$$
\begin{array}{cc}
\parallel\varrho\parallel\leq 2, & \parallel\omega\parallel\leq\frac{1}{R}
\end{array}
$$
\end{definition}

The completed tensor algebra $T_{R}A$ is natural in $A$ with respect to 
contractive algebra homomorphisms. The identity on $TA$ extends to 
a contractive homomorphism $T_RA\to T_{R'}A$ for $R<R'$. The formal inductive
limit
$$
{\cal T}A\,=\,"\underset{R\to\infty}{lim}"\,T_RA
\eqno(2.46)
$$
is called the "strict" universal infinitesimal deformation of $A$. 
Its properties are described in detail in \cite{Pu2}, section 1.

Now we can formulate

\begin{prop}
Let ${\cal E}=({\cal H},\rho,\epsilon,{\cal D})$ be an even
unbounded $\theta$-summable Fredholm module over $A$. Let
$\Phi:Q(TA)\to\widetilde{{\cal L}}$ be the characteristic morphism (2.10)
associated to $\cal E$ and let $\widetilde{\tau}$ be the odd trace on
$\widetilde{{\cal L}}$ introduced in 2.3. Then the pull-back of
$\widetilde{\tau}$ along $\Phi$ extends to a continuous trace 
$$
\Phi^*\widetilde{\tau}:Q_{R}({\cal T}A)\to\C
$$ 
on the ind-Banach algebra $Q_{R}({\cal
T}A)\,=\,"\underset{R'\to\infty}{lim}"\,Q_R(T_{R'}A)$ provided that $R>>0$
is large enough. \end{prop}
\begin{proof}
Unfortunately one cannot make use of the characteristic properties of the
involved universal algebras because for us $\widetilde{\cal L}$ is just an
abstract algebra without any distinguished topology. Therefore the continuity
of $\Phi$ does not make any sense and one can only talk about the continuity of
$\Phi^*\widetilde{\tau}$. To this end we want to make the norms on
$Q_R(T_{R'}A)$ as explicit as possible. We begin with a general remark.
If $A$ is a complex algebra and if $S\subset A$ is a subset generating $A$ as
an algebra, then there exists a largest submultiplicative seminorm on $A$
satisfying $\parallel S\parallel\leq 1$. In fact this seminorm is given by 
$$
\parallel a\parallel\,=\,\underset{a=\sum
\lambda_is_i}{Inf}\sum\vert\lambda_i\vert
\eqno(2.47)
$$
where the infimum is taken over all presentations $a=\underset{finite}{\sum}
\lambda_is_i$ with $\lambda_i\in\C$ and $s_i\in S^\infty$, the multiplicative
closure of $S$.  

Let now $A$ be a Banach algebra and fix $R\geq 1$. For $m\in\N$ let
$\parallel\,\parallel_{(R,m)}$ be the largest seminorm on $QA=A*A$ satisfying 
$$
\parallel pa^0qa^1\ldots qa^n\parallel_m\,\leq\,(2R^2+n)^m\,
R^n\parallel a^0\parallel_A\ldots\parallel a^n\parallel_A 
\eqno(2.48)
$$ 
These norms are not submultiplicative but satisfy 
$$
\parallel\alpha\alpha'\parallel_{(R,m)}\,\leq\,
\parallel\alpha\parallel_{(R,m+1)}\,\parallel\alpha'\parallel_{(R,m)}
\eqno(2.49)
$$
for all $\alpha,\alpha'\in QA$. A straightforward calculation shows that these
seminorms are related to the submultiplicative norms of 2.19
by  
$$
\parallel-\parallel_{(R,0)}\,\leq\,\parallel-\parallel_{QA_{R'}}
\eqno(2.50)
$$
provided that $R'\geq 3R^3$. In fact consider (2.47) in the case of the algebra
$QA$ and the generating subset $S=\frac{1}{R'}(\theta(U)\cup\theta^\gamma(U))$
where $U$ denotes the unit ball of $A$. Then (2.47) gives an explicit formula 
for the submultiplicative norm on $QA_{R'}$. Now for any $s=\prod s_i\in
S^\infty$ one finds $\parallel s\parallel_{(R,0)}\leq\prod\parallel
s_i\parallel_{(R,1)}\leq 1$ because $\parallel S\parallel_{(R,1)}\leq 1$.
From this the inequality
$\parallel-\parallel_{(R,0)}\leq\parallel-\parallel_{R'}$ follows.

In a similar spirit consider  for $R\geq 1$ and $m\in\N$ the tensor algebra
$TA$ equipped with the largest seminorm on $TA$ 
$\parallel-\parallel_{(R,m)}$ \cite{Pu}, 5.6, \cite{Pu2}, 1.22
which satisfies
$$
\parallel\varrho(a^0)\omega(a^1,a^2)\ldots\omega(a^{2n-1},a^{2n})\parallel_{(R,m)}\,
\leq\,(2+2n)^m\, R^{-n} \parallel
a^0\parallel_A\ldots\parallel a^{2n}\parallel_A
\eqno(2.51)
$$
Again these norms are not submultiplicative but satisfy \cite{Pu2}, 1.22
$$
\parallel\alpha\alpha'\parallel_{(R,m)}\,\leq\,
\parallel\alpha\parallel_{(R,m+1)}\,\parallel\alpha'\parallel_{(R,m)}
\eqno(2.52)
$$
One deduces as before that these seminorms on $TA$ are related to the
submultiplicative norms of (2.47) by 
$$
\parallel-\parallel_{(R,0)}\,\leq\,\parallel-\parallel_{TA_{R'}}
\eqno(2.53)
$$
provided that $R\geq 4R'$. 
The present proposition claims that for sufficiently large but
fixed $R>1$ and all $R_1>1$ the functional $\Phi^*\widetilde{\tau}$ is a
continuous trace on  $Q_R(T_{R_1}A)$. According to (2.50) this would follow
from the estimates $\vert\Phi^*\widetilde{\tau}(x)\vert\leq C(R_1)\parallel
x\parallel_{(R,0)}$. By definition of this seminorm it suffices to verify the
estimates $$
\vert\Phi^*\widetilde{\tau}(pq^k(\alpha_0,\ldots,\alpha_n))\vert\leq C(R_1)\,
R^k\,\prod\parallel\alpha_i\parallel_{T_{R_1}A}
\eqno(2.54)
$$
By (2.53) the bound (2.54) would follow from a similar estimate with respect to
the seminorms (2.51) instead of $\parallel-\parallel_{T_{R_1}A}$ so that the
claimed continuity of $\Phi^*\widetilde{\tau}$ would be finally a consequence
of the estimates  
$$
\vert\Phi^*\widetilde{\tau}(p(\varrho\omega^{n_0})q(\varrho\omega^{n_1})\ldots
q(\varrho\omega^{n_k}))\vert\leq C(R_2)\,
R^k\,\prod\parallel\varrho\omega^{n_i}\parallel_{(R_2,0)} \eqno(2.55)
$$
for $\varrho\omega^{n_0},\ldots,\varrho\omega^{n_k}\in TA$ where 
$\parallel\varrho\omega^n(a^0,\ldots,a^{2n})\parallel_{(R_2,0)}\,=\,
R_2^{-n}\prod\parallel a^i\parallel_A$. 

Recall that if
$T=T_++\lambda^{\frac12}T_-,\,\,T_\pm\in{\cal L}$ is an operator valued 
distribution in $\widetilde{\cal L}$, whose odd part possesses a majorizing
function  $T_-\underset{r}{\prec} f$, then by definition 
$$
\vert\widetilde{\tau}(\lambda^j\, T)\vert\leq C(r)^j\, j!\, f(1)
\eqno(2.56)
$$
for all $j\geq 0$.
We will make use of (2.25) and (2.26) to obtain the desired estimates. So for
$$
\alpha\,=\,p(\varrho\omega^{n_0})q(\varrho\omega^{n_1})\ldots
q(\varrho\omega^{n_k})
(a^0,\ldots,a^{k+2l}),
$$
$$
a^0,\ldots,a^{k+2l}\in
S\subset A,\,l=(n_0+\ldots+n_k)>0  
$$
and any $R\geq 1$ we obtain with $l=2r+1$ or $l=2r+2$
$$
\vert\widetilde{\tau}(\Phi(\alpha))\vert\,\leq\,(C_{34}({\cal
E})\parallel S\parallel)^{k+2l+1}\, f_{l}(1) 
$$
$$
\leq (C_{35}\parallel S\parallel)^{k+4r+1}\,
\frac{1}{r!} \,\leq\, R_2^{-2r}\, (C_{35}\parallel S\parallel)^{k+1}\,
\frac{(R_2^2(C_{35}\parallel S\parallel)^{4})^r}{r!}
$$
$$
\leq\,C_{36}({\cal
E},R_2,\parallel S\parallel)\, R_2^{-(2r+2)}\, (C_{35}({\cal
E})\parallel S\parallel)^{k+1} 
$$
If we take for $S$ the unit ball in $A$ and have a look at (2.51) and (2.55) we
see  that we just obtained the required estimate provided that $R>C_{35}({\cal
E})$. Consider now elements $\alpha\in Q(TA)$ of the form 
$$
\alpha'\,=\,p(\varrho)q(\varrho)^k(a^0,\ldots
,a^k),\,a^0,\ldots,a^{k}\in S,\,k>0  
$$
 Then according to (2.29) $\vert\Phi(\alpha')\vert\leq (C_{37}\parallel
S\parallel)^{k+1}f_1(s)$ which implies \\ $\vert\Phi(\alpha')\vert\leq 
C_{38}({\cal E})\parallel\alpha'\parallel_{(R,0)}$ for $R>C_{37}$. The 
case of the distributions $\alpha''=q(\varrho)^n$ is treated in the same way
 and the last case $\alpha=p(\varrho(a))$ follows easily from (2.29).
\end{proof}

We are able to give a first alternative description of the cohomology class of
Connes' character.
\begin{theorem} 
Let ${\cal E}\,=\,({\cal H},\rho, \epsilon, {\cal D} )$ be an even unbounded
$\theta$-summable Fredholm module over $A$ and let
$Ch_\epsilon({\cal E})\in\,CC_\epsilon^0(A)$ be its character in the sense of
Connes (see 1.5). Let $\iota_{{\cal
T}A}^\epsilon\in\,HC^\epsilon_0({\cal
T}A, q_R({\cal T}A)),\,R>>0$ be the canonical bivariant cohomology class (2.44)
associated by excision to the universal splitting extension. Let
finally $\Phi^*(\widetilde{\tau})=\Phi_*\circ\widetilde{\tau}$ be the trace
on $Q_R({\cal T}A)$ constructed in 2.21 which is well defined for $R>>0$
sufficiently large. Then   
$$
\pi_*\circ [Ch_\epsilon({\cal
E})]\,=\,\iota_{{\cal
T}A}^\epsilon\circ[\Phi_*\circ\widetilde{\tau}]\,\in\,HC^\epsilon_0({\cal
T}A,\C)\,=\,HC_\epsilon^0({\cal T}A) \eqno(2.57) 
$$ 
\end{theorem}

\begin{proof}
Because $\Phi_*\circ\widetilde{\tau}\in HC^\epsilon_0(q({\cal
T}A),\C)$ is the restriction of a trace on $Q({\cal T}A)$
one has 
$\iota_{{\cal
T}A}^\epsilon\circ\Phi_*\circ\widetilde{\tau}\,=\,(\theta_*-\theta^{\gamma}_*)\circ\Phi_*\circ\widetilde{\tau}
\,=\,\theta_*\circ\Phi_*\circ\widetilde{\tau}
-\theta^{\gamma}_*\circ\Phi_*\circ\widetilde{\tau}$.  By proposition 2.21 and
Connes' theorem 1.5 $\pi_*\circ Ch_\epsilon({\cal E})$ and
$\theta_*\circ\Phi_*\circ\widetilde{\tau}
-\theta^{\gamma}_*\circ\Phi_*\circ\widetilde{\tau}$ are entire cyclic
cocycles on ${\cal T}A$. It suffices therefore to show that they agree on the
dense subcomplex $CC_*(TA)$ of $CC_*^\epsilon({\cal T}A)$. This is true by
1.12. 
 \end{proof}

{\bf Spectra of $\theta$-summable modules}

The $\theta$-summability condition allows Connes to use a heat kernel
regularization to obtain his character formula.
We want to replace the heat kernel regularization by a cutoff regularization 
which uses spectral projections of the operator ${\cal D}^2$. Such
regularizations  occur in the (bivariant) character formulas of
\cite{Pu}, \cite{Me}. They apply to arbitrary Fredholm modules
 but do not lead to natural explicit formulas. 
In order to carry out this change of regularization some information about the
spectrum of the operator ${\cal D}^2$ is needed.

\begin{lemma}
Let ${\cal E}\,=\,({\cal H}, \rho, \epsilon,{\cal D})$ be an even unbounded 
$\theta$-summable Fredholm module. Put
$$
N_k\,:=\,\sharp\{\lambda\in Sp({\cal D}^2),\,k\leq\lambda<k+1\}
\eqno(2.58)
$$
Then for every $\alpha>0$ there exists a constant $C_\alpha$ such that 
$$
\begin{array}{cc}
N_k\,\leq\,C_\alpha\,(1+\alpha)^k, & \forall k\geq 0 \\
\end{array}
\eqno(2.59)
$$
\end{lemma}
\begin{proof}
By assumption the resolvent $(1+{\cal D}^2)^{-1}$ is a compact
selfadjoint operator. This implies that the spectrum of ${\cal D}^2$ is
discrete with finite multiplicities. Let
$0\leq\lambda_0\leq\lambda_1\leq\ldots$ be the eigenvalues  of ${\cal D}^2$
counted with multiplicities. Let $t_0>0$ be such that $e^{t_0}<1+\alpha$. Then 
$$
\sum_{k=0}^\infty N_k\,(1+\alpha)^{-k}\,\leq\,\sum_{k=0}^\infty N_k
\, e^{-t_0k} \leq\,Trace(e^{-t_0{\cal D}^2})\,<\,\infty
$$
which implies the claim.
\end{proof}
\begin{lemma}
Let ${\cal E}\,=\,({\cal H}, \rho,\epsilon,{\cal D})$ be an even unbounded
$\theta$-summable Fredholm module. There exists a sequence $P_0\leq
P_1\leq\ldots\leq P_n\leq\ldots$ of finite rank projections in ${\cal L(H)}$
which possesses the following properties.
$$
\begin{array}{cc}
\underset{n\to\infty}{\lim}\,Rank(P_n)\,=\infty, & 
\underset{n\to\infty}{\overline{\lim}}\,(Rank(P_n))^{\frac{1}{n}}\,<\,\infty
\end{array} 
$$
$$
\begin{array}{ccc}
[{\cal D},\,P_n]=0,& [\epsilon,\,P_n]=0, & \forall n\in\N
\end{array}
$$
$$
\begin{array}{c}
\underset{n\to\infty}{\lim}\,\left(\parallel(1-P_n)\, e^{-\frac{{\cal
D}^2}{4}}\parallel_1\right)^{\frac{1}{n}}\,=\,0 
\end{array}
$$
\end{lemma}

\begin{proof}
For $k\in\N,\,k>2$ consider the sequence of spectral projections \\
$P_{k,n}\,=\,P_{[0,kn[}({\cal D}^2)$. Choose $\alpha_k>0$ such that
$(1+\alpha_k)^k\leq 2$. By lemma 2.23 we find
$$
Rank(P_{k,n})\,=\,\sum_{j=0}^{kn-1}N_j\,\leq\,C_{\alpha_k}\,\sum_{j=0}^{kn-1}(1+\alpha_k)^j
\,\leq\,\alpha_k^{-1}\, C_{\alpha_k}\, 2^n
$$ 
and 
$$
Trace((1-P_{k,n})\, e^{-\frac{{\cal D}^2}{4}})\leq
\sum_{j=kn}^\infty N_j\, e^{-\frac{j}{4}}\,\leq\,C'_{\alpha_k}\,
(((1+\alpha_k)\,e^{-\frac{1}{4}})^k)^n 
$$
Choosing an appropriate subsequence of $(P_{k,n})_{k,n\in\N}$ yields a
sequence with the desired properties.
\end{proof}

We are now ready to replace the heat kernel regularization, which is central
in Connes' construction of the character of $\theta$-summable Fredholm modules 
by a more primitive cut-off regularization. This has two advantages: it
applies to Fredholm modules without any summability condition and provides much
more flexibility with respect to algebraic operations, homotopies, and so on. 
Because we want to identify Connes' character not only with the universal
bivariant Chern-Connes character in local cyclic cohomology but also with 
Meyer's Chern character in analytic cyclic cohomology, we use a change of
regularization which identifies the Connes character formula directly 
with Meyer's character and show finally that it coincides with our character 
in the local cyclic theory.

\begin{lemma}
Let ${\cal E}=({\cal H},\rho,\epsilon,{\cal D})$ be an even unbounded
$\theta$-summable Fredholm module over $A$. Choose a sequence ${\cal
P}=(P_n)_{n\in\N}$  of finite rank projections in $\cal L(H)$  which
satisfies the conditions of 2.24. Define an odd cochain 
$\Phi_{{\cal P}}^*(\widetilde{\mu}_{{\cal P}})\,=\,\Phi_{{\cal
P}}^*(\widetilde{\mu}_{{\cal P}})^{2n+1})_{n\in\N}$  on the (acyclic)
$(b,B)$-bicomplex $CC(q(TA))$ by  
$$
\Phi_{{\cal P}}^*(\widetilde{\mu}_{{\cal
P}})^{2n+1}(\alpha_0,\ldots,\alpha_{2n+1})\,=\, \widetilde{\mu}_{{\cal
P}}^{2n+1}(\Phi_{{\cal P}}(\alpha^0),\ldots,\Phi_{{\cal P}}(\alpha^{2n+1}))
\eqno(2.60) 
$$
where $\widetilde{\mu}_{{\cal P}}\in CC^1(\widetilde{{\cal L}})$ is the
cochain of 2.4 and the linear map $\Phi_{{\cal P}}:q(TA)\to\widetilde{{\cal
L}}$ was introduced in 2.17. Then for sufficiently large $R>>0$ the
functionals  $\Phi_{{\cal P}}^*(\widetilde{\mu}_{{\cal P}})\,=\,\Phi_{{\cal
P}}^*(\widetilde{\mu}_{{\cal P}})^{2n+1})_{n\in\N}$  define an entire cyclic
cochain on the ind-Banach algebra $q_R({\cal
T}A)\,=\,"\underset{R'\to\infty}{lim}"\,q_R(T_{R'}A)$.  
\end{lemma}

\begin{proof}
A reasoning similar to the one in the proof of 2.21 and the definition
\cite{Co} of entire cyclic cohomology show that it suffices to verify the
following. For some $R_0\geq1$ and all $R\geq 1$ there exist
constants $C_{2n+1}({\cal E},R)>0,\,n>0,$ such that for all bounded subsets
$S\subset A$
$$
\vert\Phi_{{\cal P}}^*(\widetilde{\mu}_{{\cal
P}})^{2n+1}(\alpha_0,\ldots,\alpha_{2n+1})\vert\,\leq\, C_{2n+1}({\cal E},R)\,
R_0^k\, R^{-l}\,\parallel S\parallel^{2n+2+2l+k} \eqno(2.61)
$$
for all elements $\alpha_0,\ldots,\alpha_{2n+1}\in q(TA)$ of the form 
$$
\alpha_i\,=\,p(\varrho\omega^{n_0})q(\varrho\omega^{n_1})\ldots
q(\varrho\omega^{n_{k_i}})(a^0,\ldots,a^{k_i+2l_i})
$$
$$
a^0,\ldots,a^{k_i+2l_i}\in
S,\,l_i=(n_0+\ldots+n_{k_i})\geq 0,\,k_i>0,
$$
$k=k_0+\ldots+k_{2n+1},\,l=l_0+\ldots+l_{2n+1}$. 
Moreover these constants have to satisfy the asymptotic decay condition
$$
\underset{n\to\infty}{\overline{lim}}\,(n!\, C_{2n+1}({\cal
E},R))^{\frac{1}{n}}\,=\,0 
\eqno(2.62)
$$
for all $R\geq 1$.
The estimates of lemma 2.18 together with 2.11 show that
$$
n!\,\vert\Phi_{{\cal P}}^*(\widetilde{\mu}_{{\cal
P}})^{2n+1}(\alpha_0,\ldots,\alpha_{2n+1})\vert\,\leq\, 
\sum_{j=0}^{n+1}\,C_{37}^j\, j!\, (C_{32}({\cal E})\parallel
S\parallel)^{k+2l+2n+2}\,\,f_{l+1}(1) 
$$
$$
\leq\,C_{38}^n\, n!\, C_{32}^k\, R^{-l}\, exp(C_{39}\,
R^2)\,\parallel S\parallel^{k+2l+2n+2} 
$$
 which gives an estimate of the type (2.61). It is however
not strong enough  to guarantee the decay condition (2.62). To achieve this we
suppose that $n>>0$ is large enough that $Im(1-P_n)\cap Ker{{\cal D}}=0$ and
$\eta_n=Trace((1-P_n)\, e^{-\frac{{\cal D}^2}{4}})\,<\,e^{-\frac{1}{4}}$. 
Then we can use the sharper estimate of (2.41) and arrive at
$$
n!\,\vert\Phi_{{\cal P}}^*(\widetilde{\mu}_{{\cal
P}})^{2n+1}(\alpha_0,\ldots,\alpha_{2n+1})\vert 
$$
$$
\leq\,
 \left(\sum_{j=0}^{n+1}\,C_{37}^j\,
j!\,(-4\,log\,\eta_n)^{-j}\right)\, (C_{33}({\cal E})\parallel
S\parallel)^{k+2l+2n+2}\,\eta_n\,\,f_{l+1}(1) 
\eqno(2.63)
$$
According to lemma 2.24 $\underset{n\to\infty}{lim}\,(\eta_n)^{\frac1n}\,=\,0$
so that $\underset{n\to\infty}{lim}\,\frac{-log\,\eta_n}{n}\,=\,\infty$. 
Thus for $n>>0$ large enough $(-4\,log\,\eta_n)\,\geq\,2C_{37}(n+1)$ and then 
$$
\sum_{j=0}^{n+1}\,C_{37}^j\,
j!\,(-4\,log\,\eta_n)^{-j}\,\leq\,\sum_{j=0}^{n+1}\,2^{-j}\,\frac{j!}{(n+1)^j}\,\leq\,2
\eqno(2.64)
$$
For the other term in (2.63) we arrive at 
$$
(C_{33}({\cal E})\parallel
S\parallel)^{k+2l+2n+2}\,\eta_n\,\,f_{l+1}(1)\,\leq\,
\left(exp(C_{40}\, R^2)\,\eta_n\right)\, R^{-l}\,
C_{33}^{2n+2+k}\,\parallel S\parallel^{k+2l+2n+2}  
$$ 
which in view of 2.24 yields the desired bound provided that $R_0\geq
C_{33}({\cal E})$. 
\end{proof}

{\bf Pointwise estimates for the operator norm of operator valued holomorphic
functions}

We take up again the study of continuity
properties of various cyclic cochains. We will not only
use uniform estimates of the Schatten norms of the involved distributions
near 0 but also uniform pointwise estimates of the operator norm of their
Laplace transforms. 

The contour $\Gamma$ of (2.4) divides the plane into two parts. Following
Connes \cite{Co} we are interested in pointwise norm estimates of various
operator valued holomorphic functions on the closure of the
connected component of $\C\backslash\Gamma$ contained in the right halfplane.
Explicitely it is given by 
$$
\overline{U}\,=\,\{z\in\C,\,\vert
Im(z)\vert\,\leq\,\sqrt{3}\,Re(z-1)\} \eqno(2.65)
$$

\begin{lemma}
Let ${\cal E}=({\cal H}, \rho , \epsilon, {\cal D})$ be an even unbounded 
$\theta$-summable Fredholm module over $A$ and let $F(z)\,=\,({\cal
D}\,+\,z^{\frac12}\,\epsilon)({\cal D}^2\,+\,z)^{-\frac12}$ be the Laplace
transform of the associated involutive distribution. Then
$$
\begin{array}{cc}
\parallel F(z)\parallel_{{\cal L(H)}}\,\leq\,C_{41}, &
\parallel [F(z),a]\parallel_{{\cal
L(H)}}\,\leq\,C_{42}({\cal E})\,\vert z\vert^{-\frac12}\parallel a\parallel_A
\end{array}
$$
for all $a\in A$ and $z\in\overline{U}^2$. Similar estimates hold for the
Laplace transforms of the even respectively odd parts (with respect to the
grading of $\widetilde{{\cal L}}$)  of the involved distributions.
\end{lemma}

\begin{proof}
The decomposition of the operator under consideration into even and odd parts
$$
\begin{array}{ccccc}
F(z) & = & ({\cal D}\,+\,\epsilon z^{\frac12})({\cal D}^2\,+\,z)^{-\frac12}
& = &  F_1(z^{-\frac12}{\cal D})\,+\epsilon\,F_2(z^{-\frac12}{\cal D}) \\
\end{array}
\eqno(2.66)
$$
is given by operators obtained from ${\cal D}'=z^{-\frac12}{\cal D}$ by
holomorphic functional calculus with respect to
$$
\begin{array}{cc}
F_1(u)\,=\,u\,(1+u^2)^{-\frac12}, &
F_2(u)\,=\,(1+u^2)^{-\frac12}\\
\end{array}
\eqno(2.67)
$$
By assumption the operator $\cal D$ is selfadjoint
and its resolvent $(1+{\cal D}^2)^{-1}\in{\cal K(H)}$ is compact. Therefore
the spectrum of $\cal D$ is real and discrete with finite multiplicities. 
One deduces easily that there are universal
bounds $$
\begin{array}{cc}
\parallel F_1(z^{-\frac12}{\cal D})\parallel_{{\cal L(H)}}\,\leq\,C_{43}, & 
\parallel F_2(z^{-\frac12}{\cal D})\parallel_{{\cal L(H)}}\,\leq\,C_{44}
\end{array}
\eqno(2.68)
$$
for $z\in\overline{U}^2$ from which the first assertion follows.
Our analysis of the commutators $[F(z),a],\,a\in A,$ is based on the identity
\cite{BJ} 
$$
\begin{array}{ccc}
(1\,+\,{\cal D}^2)^{-\frac12} & = &
\frac{1}{\pi}\int_0^\infty\lambda^{-\frac12}(1\,+\,{\cal
D}^2\,+\,\lambda)^{-1}\,d\lambda \\
\end{array}
\eqno(2.69)
$$
For the commutator of the even part of $F$ and $a\in A$ we find
$$
[F(z)_+,a]\,=\,[F_1(z^{-\frac12}{\cal D}), a]\,=
\,\frac{1}{\pi}\int_0^\infty\lambda^{-\frac12}
[{\cal D'},a](1+{\cal D'}^2+\lambda)^{-1}\,d\lambda
$$
$$
-\frac{1}{\pi}\int_0^\infty\lambda^{-\frac12}
({\cal D'}(1+{\cal D'}^2+\lambda)^{-1}[{\cal D'},a]
{\cal D'}(1+{\cal D'}^2+\lambda)^{-1}\,d\lambda
$$
$$
-\frac{1}{\pi}\int_0^\infty\lambda^{-\frac12}
{\cal D'}^2(1+{\cal D'}^2+\lambda)^{-1}[{\cal D'},a]
(1+{\cal D'}^2+\lambda)^{-1}\,d\lambda
$$
To estimate the norms of these operators we observe that for 
all $u\in \overline{U}$
and all $\lambda\in\R_+$ one has 
$\vert (1+u^2+\lambda)^{-1}\vert\,\leq\,2\,(1+\lambda)^{-1},\,
\vert u(1+u^2+\lambda)^{-1}\vert\,\leq\,2\,(1+\lambda)^{-\frac12}
$ and $\vert u^2(1+u^2+\lambda)^{-1}\vert\,\leq\,2$. With the help of these 
bounds we deduce 
$$
\parallel[F(z)_+,a]\parallel\,\leq
\,C_{45}\parallel[{\cal D'},a]\parallel\,=\,
C_{46}\vert z\vert^{-\frac12}\parallel[{\cal D},a]\parallel
\eqno(2.70)
$$
for all $z\in\overline{U}^2$. The corresponding assertion 
for the odd part of $F(z)$ follows from similar arguments.
\end{proof}

\begin{lemma}
Let ${\cal E}=({\cal H}, \rho, \epsilon,{\cal D})$ be an even unbounded 
$\theta$-summable Fredholm module over $A$ and let
$\Phi:q(TA)\to\widetilde{{\cal L}}$ be its characteristic homomorphism
introduced in 2.6. Let $P$ be a spectral projection associated to ${\cal
D}^2$ and suppose  that $Im(P)\cap ker({\cal D})=0$ and that
$\eta\,=\,Trace(P\, e^{-\frac{{\cal D}^2}{4}})\,<\,1$. Then the image  of the
distribution $\Phi(q(\varrho(a)) P,\,a\in A,$ under the Laplace transform ,as
well as its even and odd parts satisfy  
$$
\parallel\Phi(q(\varrho(a))(z) P\parallel\,\leq\,C_{47}\left(
(-log\,\eta)^{-\frac12}\vert z^{\frac12}\vert\,+\,\vert
z\vert^{-\frac12}\right)\parallel a\parallel \eqno(2.71) 
$$
for all $a\in A,\,z\in\overline{U}^2$.
\end{lemma}
\begin{proof}
We calculated in (2.30) the distribution $\Phi(q(\varrho(a))$. In fact 
$$
\Phi(q(\varrho(a))(z)\,
P\,=\,\left(\frac12\,a\,F(z)_-\,+\,\frac14[F(z)_-,\,a]\,-\,
\frac14F(z)_-[{\cal D},\,a]({\cal D}^2+z)^{-\frac12} \right) P 
$$ 
According to 2.26 $\parallel \frac14[F(z),\,a]_-\,-\, \frac14F(z)_-[{\cal
D},\,a]({\cal D}^2+z)^{-\frac12}
\parallel\,\leq\,C_{48}\vert z\vert^{-\frac12}\parallel a\parallel$. On the
other hand $\parallel F(z)_-P\parallel\,=\,\parallel z^{\frac12}\epsilon({\cal
D}^2+z)^{-\frac12}P\parallel\,\leq\,2\vert
z\vert^{\frac12}\,\lambda_1^{-\frac12}$ where $\lambda_1$ is the smallest
nonzero eigenvalue of ${\cal D}^2P$.  It satisfies $\lambda_1\geq
-4\,log\,\eta$ (see the proof of 2.12) whence the result. 
\end{proof}

The estimates obtained so far yield

\begin{prop}
Let ${\cal E}=({\cal H},\rho,\epsilon,{\cal D})$ be an even unbounded
$\theta$-summable Fredholm module over $A$ and let
$\Phi:q(TA)\to\widetilde{{\cal L}}$ be  its characteristic homomorphism (see
(2.10)). Choose a sequence of finite rank projections in $\cal L(H)$  which
satisfies the conditions of 2.24 and consider the odd cochain
$\widetilde{\mu}_{{\cal P}}\,=\,(\mu^{2n+1})_{n\in\N}\in CC^1$ (2.6), 
on the $(b,B)$-bicomplex of $\widetilde{{\cal L}}$. 
Then the pull-back of this functional along $\Phi$ extends to an entire cyclic
cochain 
$$
\Phi^*(\widetilde{\mu}_{{\cal P}})\,\in\,CC_\epsilon^1(q_{R_0}({\cal T}A))
$$
on the ind-Banach algebra $q_{R_0}({\cal
T}A)\,=\,"\underset{R\to\infty}{lim}"\,q_{R_0}(T_RA)$ provided that
$R_0>>0$ is large enough.
\end{prop}

\begin{proof}
It will not be shown directly that $\Phi^*(\widetilde{\mu}_{{\cal P}})$ is an entire
cyclic cochain  but we will rather prove that the difference
$\Phi^*(\widetilde{\mu}_{{\cal P}})\,-\,\Phi_{{\cal
P}}^*(\widetilde{\mu}_{{\cal P}})$ (see 2.17) is entire. The conclusion
follows then in conjunction with 2.25. As we have seen in the proofs of 2.21
and 2.25 it suffices to verify the following. There exists some $R_0\geq 1$
such that for all $R_1,\,R_2\geq 1$ there are constants $C({\cal
E},R_1,\,R_2)>0$ so that for every bounded subset $S\subset A$  
$$
n!\,\vert(\Phi^*(\widetilde{\mu}_{{\cal P}})^{2n+1}\,-\,\Phi_{{\cal
P}}^*(\widetilde{\mu}_{{\cal
P}})^{2n+1})(\alpha_0,\ldots,\alpha_{2n+1})\vert\,\leq\, C(R_1,R_2)\, R_0^k\,
R_1^{-l}\, R_2^{-n}\, \parallel S\parallel^{2n+2+2l+k} \eqno(2.72)
$$
for all $n\geq 0$ and all elements $\alpha_0,\ldots,\alpha_{2n+1}\in q(TA)$ of
the form  
$$
\alpha_i\,=\,p(\varrho\omega^{n_0})q(\varrho\omega^{n_1})\ldots
q(\varrho\omega^{n_{k_i}})(a^0,\ldots,a^{k_i+2l_i})
$$
$$
a^0,\ldots,a^{k_i+2l_i}\in
S,\,l_i=(n_0+\ldots+n_{k_i})\geq 0,\,k_i>0,
$$
$k=k_0+\ldots+k_{2n+1},\,l=l_0+\ldots+l_{2n+1}$. 

To understand the algebraic manipulations which follow note that for any 
idempotent $P\in A$ and any elements $a^0,\ldots,a^n\in A$ one has the
identity
$$
P(\prod_0^n a^i)P\,-\,\prod_0^n(Pa^iP)\,=\,
\sum_{j=0}^{n-1}\,P(\prod_0^j
a_i)[a^{j+1},P](\prod_{j+2}^nPa_kP)P
\eqno(2.73)
$$
In the sequel the notations of lemma 2.17 are understood. 
For \\ $\alpha\,=\,p(\varrho\omega^{n_0})q(\varrho\omega^{n_1})\ldots
q(\varrho\omega^{n_{k}})(a^0,\ldots,a^{k+2l}),\,
a^0,\ldots,a^{k+2l}\in
S,$ \\ $l_i=(n_0+\ldots+n_{k})\geq 0,\,k>0$ one finds 
$$
n!\,\,\Phi^*(\widetilde{\mu}_{{\cal
P}})^{2n+1}(\alpha_0,\ldots,\alpha_{2n+1})\,=\,
\widetilde{\tau}\left(\prod_{i=0}^{2n+1}
(\Phi^{(1-P_n)}(\alpha_i)-\beta_i)\right) \eqno(2.74)
$$
$$
n!\,\,\Phi_{{\cal
P}}^*(\widetilde{\mu}_{{\cal
P}})^{2n+1}(\alpha_0,\ldots,\alpha_{2n+1})\,
=\,\widetilde{\tau}\left(\prod_{i=0}^{2n+1} \Phi^{(1-P_n)}(\alpha_i)\right) 
$$ 
$$
\beta\,=\,\Phi^{(1-P_n)}(\alpha)\,-\,(1-P_n)\,\Phi(\alpha)\,(1-P_n)\,=
$$
$$
=\,(1-P_n)\Phi'(p)
\Phi^{(1-P_n)}(q^{k})(1-P_n)\,+
$$
$$
+\,(1-P_n)\sum_{j=1}^k\Phi(pq^{j-1})\Phi'(q)
\Phi^{(1-P_n)}(q^{k-j})(1-P_n)
\eqno(2.75)
$$
$$
\begin{array}{cc}
\Phi'(p)\,=\,\frac12(T\varphi_0''+T\varphi_1''),
& \Phi'(q)\,=\,\frac12(T\varphi_0''-T\varphi_1'') \\
\end{array}
$$
$$
T\varphi_0''(\varrho\omega^m)(a^0,\ldots,a^{2m})\,=\,
C_n\,[P_+a^0P_+,\,Q_n]\prod_1^{2m}(Q_n[F,a^i]Q_n)\,+
$$
$$
+\,C_n\sum_{m'=0}^{2m-1}(P_+a^0P_+)\,\left(\prod_1^{m'-1}[F,a^{i}]\right)\,
[[F,a^{m'}],Q_n]\left(\prod_{m'+1}^{2m}(Q_n[F,a^{i}]Q_n)\right)
\eqno(2.76)
$$
$$
T\varphi_1''(\varrho\omega^m)(a^0,\ldots,a^{2m})\,=\,
C_n\,\epsilon[
P_-a^0P_-,\,Q_n]\prod_1^{2m}(Q_n[F,a^i]Q_n)\epsilon\,+
$$
$$
+\,C_n\sum_{m'=0}^{2m-1}\epsilon(P_-a^0P_-)\,\left(\prod_1^{m'-1}[F,a^{i}]\right)\,
[[F,a^{m'}],Q_n]\left(\prod_{m'+1}^{2m}(Q_n[F,a^{i}]Q_n)\right)\epsilon
\eqno(2.77)
$$
for $m>0$ where $Q_n\,=\,(1-P_n)$ and $C_n=\left(-\frac{1}{4}\right)^n$. For
$m=0$  $$
\begin{array}{cc}
T\varphi_0''(\varrho)(a^0)\,=\,[P_+a^0P_+,\,(1-P_n)], & 
T\varphi_1''(\varrho)(a^0)\,=\,[\epsilon P_-a^0P_-\epsilon,\,(1-P_n)] \\
\end{array}
\eqno(2.78)
$$
The fact that all commutators with the projection $(1-P_n)$ are of finite rank 
at most $2\,Rank(P_n)$ enables us to use the Laplace transform and its inverse 
to control the trace of distributions of the form 
$T\,=\,T_0[T_1,(1-P_n)]T_2,\,T_i\in\widetilde{{\cal L}}$. Let
$T_i=T_{i+}+T_{i-}$ be the decomposition of $T_i$ into
even and odd parts. Then one finds for any oriented curve $\Gamma'$ in
$\overline{U}$ (2.65), which is homologous to $\Gamma=\partial\overline{U}$,
the estimate  
$$
\vert\widetilde{\tau}(T)\vert\,\leq\,\sum\,2\,Rank(P_n)\,\underset{\Gamma'}{\int}\,
\parallel T_{0\pm}(z^2)\parallel\, \parallel
T_{1\pm}(z^2)\parallel\, \parallel
T_{2\pm}(z^2)\parallel\,\vert d\mu(z)\vert \eqno(2.79)
$$
where the sum runs over all possible choices of even resp. odd parts of the
involved operator valued functions with odd total weight. Here
$\parallel-\parallel$ denotes the pointwise operator norm of an operator
valued holomorphic function. 

If we apply this argument in our case with the curve $\Gamma_\zeta$ given by
the oriented boundary of $U\cap\{z\in\C,\,Re(z)>\zeta\}$ we find the following
estimate. 
$$
n!\,\,\vert(\Phi^*(\widetilde{\mu}_{{\cal
P}})^{2n+1})\,-\,\Phi_{{\cal
P}}^*(\widetilde{\mu}_{{\cal
P}})^{2n+1})(\alpha_0,\ldots,\alpha_{2n+1})\vert\,\leq
$$
$$
\leq\,2\,Rank(P_n)\,
(C_{49}({\cal E})\parallel
S\parallel)^{k+2l+2n+2}\,\sum_{i,j}\zeta^{-(2l+i)}\,
(-log\,\eta_n)^{-\frac{j}{2}} \underset{\Gamma_\zeta}{\int}\,\vert
z^2\vert^{\frac{j}{2}}\vert d\mu(z)\vert 
\eqno(2.80) 
$$
where the sum runs over the integers $i,j\geq 0,0\leq i+j\leq 2n+2,\,i+j+l\geq
2n+2$. Put $R_0=C_{49}({\cal E})$. If one chooses $\zeta$ such that
$\zeta>C_{49}\, R_1$ one sees  that for each $n\geq0$ there exist
constants $C(n,R_1,R_2)$ such that the estimate (2.72) holds. What remains to
be done is to show that there is a similar bound which is independent of $n$.
It suffices to obtain uniform estimates  for $n>>0$. Note that by construction
(see 2.24) $Rank\,\, P_n\,\leq\,C_{50}^n$ for a suitable constant $C_{50}$.
Choose now $\zeta>R_1\, R_2\, C_{49}^2\, C_{50}$. As
$\underset{n\to\infty}{lim} \eta_n^{\frac1n}\,=\,0$ we may suppose that for
$n>>0$ large enough $-\,log\,\eta_n \geq\,4\,C_{49}^2\, C_{50}^2\,
R_2\,n$. Then we deduce from (2.80) 
$$
n!\,\,\vert(\Phi^*(\widetilde{\mu}_{{\cal
P}})^{2n+1}\,-\,\Phi_{{\cal
P}}^*(\widetilde{\mu})^{2n+1})(\alpha_0,\ldots,\alpha_{2n+1})\vert\,\leq
$$
$$
\leq\,2\,\, R_0^k\, R_1^{-l}\, R_2^{-n}\parallel
S\parallel^{k+2l+2n+2}\,\sum_{i,j}\underset{\Gamma_\zeta}{\int}\,\left(\frac{\vert
z^2\vert}{4n}\right)^{\frac{j}{2}}\vert d\mu(z)\vert 
\eqno(2.81) 
$$
If one can show that 
$$
\underset{n>0}{sup}\,\,\,\underset{0\leq j\leq
2n+2}{max}\,\,\,\underset{\Gamma_\zeta}{\int}\,\left(\frac{\vert
z^2\vert}{4n}\right)^{\frac{j}{2}}\vert d\mu(z)\vert\,<\,\infty 
\eqno(2.82)
$$
for all $\zeta\geq 1$ then the desired estimate (2.72) follows. 

We may suppose $\frac{\vert
z^2\vert}{4n}>1$. Then (remember that $j\leq 2n+2$)
one has  $(\frac{\vert
z^2\vert}{4n})^{\frac{j}{2}}\,\leq\,(\frac{\vert
z^2\vert}{4n})^{n+1}\,\leq\,\frac{\vert
z^2\vert^{n+1}}{4^{n+1}(n+1)!}\,\frac{(n+1)!}{n^{n+1}}\,\leq\,2\,e^{\frac{\vert
z^2\vert}{4}}$. 
It is easily seen that $Re(z^2)\,\leq\,-\frac{\vert
z^2\vert}{3}$ for $z\in\Gamma_\zeta$ provided that $\vert
z^2\vert\,>>0$ is large. Thus finally   
$$
\underset{n>0}{sup}\,\,\,\underset{0\leq j\leq
2n+2}{max}\,\,\,\underset{\Gamma_\zeta}{\int}\,\left(\frac{\vert
z^2\vert}{4n}\right)^{\frac{j}{2}}\vert d\mu(z)\vert\,
\leq\,\underset{\Gamma_\zeta}{\int}\,e^{\frac{\vert z^2\vert}{4}}\,\vert
d\mu(z)\vert
$$
$$
\leq\,\frac{1}{\pi}\,\underset{\Gamma_\zeta}{\int}\,e^{(\frac{\vert
z^2\vert}{4}+Re(z^2))} \vert
dz\vert\,\leq\,\frac{1}{\pi}\,
\underset{\Gamma_\zeta}{\int}\,e^{-\frac{\vert
z^2\vert}{12}} 
\vert dz\vert\,+\,O(1)\,<\,\infty
$$
which proves the proposition.
\end{proof}

\begin{theorem} 

Let ${\cal E}\,=\,({\cal H},\rho , \epsilon, {\cal D})$ be an even unbounded
$\theta$-summable Fredholm module over $A$ and let $[Ch_\epsilon({\cal E})]\in
HC_\epsilon^0(A)$ be the class of its Connes character 1.5 in entire cyclic
cohomology. Choose a sequence $\cal P$ of finite rank projections
in  $\cal L(H)$ satisfying the conditions of 2.24 and let
$\widetilde{\tau}_{{\cal P}}\in CC^0(\widetilde{{\cal L}})$ be the
associated cocycle (2.6).  
Let $\Phi:q(TA)\to\widetilde{{\cal L}}$ be the characteristic homomorphism
(2.10) associated to ${\cal E}$ and denote finally by $\iota_{{\cal
T}A}^\epsilon\,\in\,HC^\epsilon_0({\cal T}A, q_R({\cal T}A)),\,R>>0$ the
canonical bivariant cohomology class (2.44) associated to the universal
splitting extension. Then for $R>>0$ sufficiently large the functional
$\Phi^*(\widetilde{\tau}_{{\cal P}})$ extends to an entire cyclic cocycle on
the ind-Banach algebra $q_R({\cal T}A)$ and its cohomology class
$[\Phi_*\circ\widetilde{\tau}_{{\cal P}}]\in HC^\epsilon_0(q_R({\cal T}A),\C)$
satisfies     $$
\pi_*\circ[Ch_\epsilon({\cal E})]\,=\,\iota_{{\cal
T}A}^\epsilon\circ[\Phi_*\circ\widetilde{\tau}_{{\cal P}}]\,\in\,HC^\epsilon_0({\cal
T}A,\C)\,=\,HC_\epsilon^0({\cal T}A) \eqno(2.83)  
$$ 
\end{theorem}

\begin{proof}
This is an immediate consequence of 2.22 and 2.28.
\end{proof}

\section{ The basic conjugation and homotopy}

In the previous section we obtained a description of Connes' character of an
unbounded $\theta$-summable Fredholm module 1.5, which is formally similar to
that of the  abstract character in $K$-homology in (1.27). There is 
a basic difference however: if one compares the characteristic homomorphisms
(1.47) and (1.32) used to define the characters in 2.29 and (1.39), then one
observes that the roles of the operators $F$ and $\epsilon$ have been exchanged
in the corresponding formulas. Modulo compact operators the
data defining a bounded Fredholm module are in fact completely symmetric in $F$
and $\epsilon$. A bounded Fredholm module yields a representation of
the Clifford algebra  of the two-dimensional vector space with orthonormal
basis $\{F,\epsilon\}$ on the Calkin algebra $\cal L(H)/K(H)$. There
exists an inner automorphism  of the Clifford algebra which interchanges the
basis elements $F$ and $\epsilon$. It is given by  conjugation with the
involutive element $U\,=\,\frac{1}{\sqrt{2}}(F+\epsilon)$. We will study this
element in the algebra $\cal O$ of holomorphic operator valued functions on
$\C\backslash\R_-$ and will show that conjugation by  this element
interchanges the characteristic homomorphisms (1.47) and (1.32) modulo a small
error. As traces are invariant under conjugation we will obtain a formula for
Connes' character which coincides with that of the abstract character
modulo a small error. This error can be eliminated by a straightforward
homotopy and we arrive at our final description 3.8 of Connes' character
class in entire cyclic cohomology.

\begin{definition}
Let ${\cal E}\,=\,({\cal H},\,\rho,\,\epsilon,\,{\cal D})$ be an even unbounded
$\theta$-summable Fredholm module over $A$
and let $F(z),\,Re(z)>0,$ be the associated holomorphic family 
of operators (1.5). We define a characteristic family of operators by
$$
U(z)\,=\,\frac{1}{\sqrt{2}}(F(z)+\epsilon)
\eqno(3.1)
$$
It satisfies
$$
\begin{array}{cc}
U(z)\, F(z)\,=\,\epsilon\, U(z), & U(z)\,\epsilon\,=\,F(z)\, U(z)\\  
\end{array}
\eqno(3.2)
$$
\end{definition}

\begin{remark}
There is an element $U=\frac{1}{\sqrt{2}}(F+\epsilon)\in\widetilde{{\cal L}}$
in the associated  algebra of operator valued distributions. It is 
invertible in $\widetilde{{\cal L}}$ if and only if $Ker\,{\cal D}\,=\,0$.
Therefore conjugation with this element
makes no sense among operator valued distributions with controlled
growth if $\cal D$ possesses a kernel. It turns out however that its Laplace
transform $U(z)$ is invertible among bounded operator valued  holomorphic
functions (note that the Laplace transform is one to one on $\cal L$ but not on
$\widetilde{{\cal L}}$).
\end{remark}

\begin{lemma}
Let ${\cal E}\,=\,({\cal H},\,\rho,\,\epsilon,\,{\cal D})$ be an even unbounded
$\theta$-summable Fredholm module and let $U(z)$ be the
characteristic family of operators (3.1). Then $U(z)$ is invertible for all
$z\in\C\backslash\R_-$ and  
$$
\begin{array}{cc}
\parallel U(z)\parallel_{{\cal L(H)}}\,\leq\,C_{51}({\cal E}), & 
\parallel U(z)^{-1}\parallel_{{\cal L(H)}}\,\leq\,C_{52}({\cal E}),
 \\ \end{array}
\eqno(3.3)
$$
$$
\parallel [U(z),a]\parallel_{{\cal L(H)}}\,\leq\,C_{53}({\cal E})\,\vert
z\vert^{-\frac12}\parallel a\parallel_A
\eqno(3.4)
$$
for all $z\in \overline{U}^2$ and $a\in A$.
\end{lemma}

\begin{proof}
Observe that $[U(z),a]\,=\,\frac{1}{\sqrt{2}}[F(z),a]$
and 
$$
U^2(z)\,=\,\frac12(2\,+\,(\epsilon\,F(z)\,+\,F(z)\,\epsilon))\,=\,
1\,+\,z^{\frac12}(z+{\cal D}^2)^{-\frac12}\,=\,1+F_2(z^{-\frac12}{\cal D})
$$
in the notations of 2.26. The proof is then similar to the one of 2.26. 
\end{proof}

\begin{definition}
Let ${\cal E}\,=\,({\cal H},\,\rho,\,\epsilon,\,{\cal D})$ be an even
unbounded $\theta$-summable Fredholm module over $A$ and let $F(z)$
and  $U(z)$ be the associated families of operators (1.5). Let $\cal O$ be the
algebra of holomorphic operator valued functions on $\C\backslash\R_-$
 and define  families of linear maps 
$\rho^t,\rho_{0,1}^t:A\to{\cal O}[t]$ by 
$$
\rho^t(a)\,=\,\rho(a)\,+\,t\,[U(z),\rho(a)]\, U(z)^{-1}
\eqno(3.5)
$$
$$
\begin{array}{cc} 
\rho_0^t(a)\,=\,\frac{1+\epsilon}{2}\,\rho^t(a)\,\frac{1+\epsilon}{2}, & 
\rho_1^t(a)\,=\,F(z)\,\frac{1-\epsilon}{2}\,\rho^t(a)\,\frac{1-\epsilon}{2}\,
F(z) \\ 
\end{array}
\eqno(3.6)
$$
The maps $\rho_0^t$ and $\rho_1^t$ coincide modulo the ideal $\cal
O_K\subset\cal O$ of holomorphic functions 
with values in the ideal of compact operators.
\end{definition}

One has the 

\begin{prop}
Let ${\cal E}\,=\,({\cal H},\,\rho,\,\epsilon,\,{\cal D})$ be an even unbounded
$\theta$-summable Fredholm module and let ${\cal
P}=(P_n)_{n\in\N}$ be a family of finite rank projections satisfying the
conditions of 2.24. For $z\in\C\backslash\R_-$ let $\vartheta_{{\cal
P}}(z)\in CC^0({\cal O_K})$ be the pullback of the cocycle $\tau_{{\cal
P}}\in CC^0({\cal K(H)})$ (1.37) along the evaluation homomorphism at $z$.
Let 
$$
\Psi_t\,=\,T\rho_0^t*T\rho_1^t:\,q(TA)\,\to\,{\cal O_K}[t]
\eqno(3.7)
$$
be the family of homomorphisms of universal algebras associated to the
linear maps of 3.4. Then for fixed $R>>0$ the
following holds: for each $t\in [0,1]$ the cochain 
$$
\nu_t\,=\,
\underset{\Gamma}{\int}\Psi_t^*(\tau_{\cal
P}(z^2))\,d\mu(z)\in CC^0(q(TA))
\eqno(3.8)
$$ 
extends to an entire cyclic cocycle 
$$
\nu_t\in CC_\epsilon^0\,(q_R({\cal T}A))
\eqno(3.9)
$$
The cohomology class of this family of cocycles is independent of $t$.
\end{prop}

\begin{proof}
The map $\rho^t$ is the sum of the structure homomorphism $\rho$ and a 
linear perturbation term $t\,[U(z),\rho(a)]\,
U(z)^{-1}=\frac{t}{\sqrt{2}}\,[F(z),\rho(a)]\, U(z)^{-1}$ whose
pointwise  operator-norm is bounded according to 3.3 by 
$\parallel [U(z),\rho(a)]\,
U(z)^{-1}\parallel\,\leq\,C_{53}\,\vert z\vert^{-\frac12}$. As one may replace
the  path of integration $\Gamma$ by any of the paths $\Gamma_\zeta$ of
(2.80), one may suppose that the perturbation terms $\rho^t-\rho$ are
arbitrarily small  in norm. For $t=0$ a calculation similar to the one at the
end of the proof of 2.28 allows to deduce that  $\nu_0$ is
an entire cyclic cocycle. The case $t\neq 0$ follows  by a calculation 
making use of the case $t=0$  and the remark about the norm
of the linear perturbation term. The last  assertion is a consequence of the
Cartan homotopy formula in entire cyclic cohomology.
\end{proof}

We are going to identify the cohomology class constructed in the previous
proposition.

\begin{lemma}
Let ${\cal E}\,=\,({\cal H},\,\rho,\,\epsilon,\,{\cal D})$ be an even unbounded
$\theta$-summable Fredholm module over $A$. Let
$\Phi^*(\widetilde{\tau}_{{\cal P}})$ and $\nu_t,t\in[0,1],$ be the entire
cyclic cocycles over $q_R({\cal T}A)$  introduced in 2.29 and 3.5.,
respectively. Then $\Phi^*(\widetilde{\tau}_{{\cal P}})=\nu_1$ as entire
cyclic cocycles.  
\end{lemma}

\begin{proof}
It suffices to check that both cocycles coincide on the dense subcomplex
$CC_0(q(TA))$ of $CC^\epsilon_0\,(q_R({\cal T}A))$. After applying the
Laplace transform one sees that the first cocycle under discussion can be
written as
$$
\Phi^*(\widetilde{\tau}_{{\cal
P}})\,=\,\underset{\Gamma}{\int}\Phi^*(\tau_{\cal P}(z^2))\,d\mu(z)
$$
where $\Phi:q(TA)\to{\cal O_K}$ is the
characteristic homomorphism (1.48) of $\cal E$. Denote by $c(U(z))$ the
automorphism of $\cal O_K$ given by conjugation  with the characteristic
operator valued function $U(z)$ of (3.1). Because the operator $U(z)$
commutes with ${\cal D}^2$ it commutes  with the projections $(P_n)_{n\in\N}$
of $\cal P$. The cocycle $\tau_{{\cal P}}$ is thus invariant under
conjugation by $U(z)$: $c(U(z))^*\tau_{{\cal P}}=\tau_{{\cal P}}$.
So we find $\Phi_* \circ \tau_{\cal
P}(z)\,=\,\Phi_* \circ c(U(z))_* \circ \tau_{{\cal P}}(z)
\,=\,(c(U(z)) \circ \Phi)_* \circ \tau_{\cal P}(z)$. Now
$c(U(z))\circ\Phi\,=\,(c(U(z))\circ T\varphi_0)* (c(U(z))\circ
T\varphi_1)\,=\,T(c(U(z))\circ \varphi_0)*T(c(U(z))\circ \varphi_1)$ because
$c(U(z))$ is a homomorphism of algebras. (The notations are those of (1.47) and
(1.32)). Further
$$
(c(U(z))\circ\varphi_0)(a)\,=\,U(z)\,\frac{1+F(z)}{2}\,\rho(a)\,\frac{1+F(z)}{2}
\, U(z)^{-1}\,=
$$
$$
=\,\frac{1+\epsilon}{2}\, U(z)\,\rho(a)\,
U(z)^{-1}\frac{1+\epsilon}{2}\,=\,(\rho_0\circ c(U(z)))(a)
$$
and similarly $c(U(z))\circ\varphi_1\,=\,\rho_1\circ c(U(z))$. Finally
$$
U(z)\rho(a)U(z)^{-1}\,=\,\rho(a)\,+\,[U(z),\rho(a)]U(z)^{-1}\,=\,\rho^1(a)
$$
so that
$$
\begin{array}{cc}
c(U(z))\circ\varphi_0\,=\,\rho_0^1, & c(U(z))\circ\varphi_1\,=\,\rho_1^1
\end{array}
$$
in the notations of 3.4. 
Altogether we have shown that
$c(U(z))\circ\Phi\,=\,\Psi_1:\,q(TA)\,\to\,{\cal O_K}$. Therefore
$$
\Phi^*(\tau_{\cal
P}(z))\,=\,(c(U(z))\circ\Phi)^*(\tau_{\cal
P}(z))\,=\,(\Psi_1)^*(\tau_{\cal P}(z))
$$
and the lemma is proved.
\end{proof}

\begin{lemma}
Let ${\cal E}\,=\,({\cal H},\,\rho,\,\epsilon,\,{\cal D})$ be an even unbounded
$\theta$-summable Fredholm module over $A$ and let
$F(z),\,z\in\C\backslash\R_-$ be the associated holomorphic family of
operators. Let  
$$
\begin{array}{ccc}
\rho_0,\rho_1:\,A\to {\cal O}, &
\rho_0(a)\,=\,\frac{1+\epsilon}{2}\,\rho(a)\,\frac{1+\epsilon}{2}, & 
\rho_1(a)\,=\,F(z)\,\frac{1-\epsilon}{2}\,\rho(a)\,\frac{1-\epsilon}{2}\,
F(z) \\ \end{array}
$$
be the characteristic homomorphisms of the corresponding family of bounded 
Fredholm modules ${\cal E}(z)$ (1.5) and denote by
$\psi\,=\,\rho_0*\rho_1:\,qA\to{\cal O}$ their free product. Let finally
$\pi:\,{\cal T}A\to A$ be the canonical epimorphism. Then for $R>>0$
sufficiently large the cochain 
$$
\beta\,=\,
\underset{\Gamma}{\int}\psi^*(\tau_{\cal
P}(z^2))\,d\mu(z)\in CC^0(q(A))
\eqno(3.10)
$$ 
extends to an entire cyclic cocycle 
$$
\beta\in CC_\epsilon^0\,(q_R(A))
\eqno(3.11)
$$
which satisfies 
$$
q_R(\pi)_*\,\circ\,\beta\,=\,\nu_0\,\in\,CC_\epsilon^0(q_R({\cal T}A))
\eqno(3.12)
$$
\end{lemma}

\begin{proof}
This is a consequence of the fact that the linear maps $\rho_{0,1}^t:A\to{\cal
O}$ of 3.4 are actually homomorphisms of algebras for $t=0$. Therefore the
characteristic map factors as
$\Psi_0:q(TA)\overset{q(\pi)}{\longrightarrow}qA\overset{\psi}
{\longrightarrow}{\cal O_K}$. The estimates used already in the proof of 3.5
yield then the desired result.  
\end{proof}

We arrive now at the final formula for an entire cyclic cocycle representing 
the cohomology class of Connes' character 1.5. The original unbounded 
Fredholm module enters only through the characters of its associated
holomorphic family  of bounded Fredholm modules. This makes it very easy to
compare Connes' character with various Chern-Connes characters of
bounded Fredholm modules.

\begin{theorem}
Let ${\cal E}\,=\,({\cal H},\,\rho,\,\epsilon,\,{\cal D})$ be an even 
unbounded $\theta$-summable Fredholm module over $A$. Let ${\cal
P}=(P_n)_{n\in\N}$  be a family of finite rank projections satisfying the
conditions of 2.24 and let $\tau_{\cal P}$ be the corresponding analytic
cyclic cocycle over $\cal K(H)$ (see 1.37).
\begin{itemize}
\item[a)] For $R>>0$ sufficiently large denote by 
$\iota^\epsilon_A:\,CC_*^{\epsilon}(A)\,\to\,
CC_*^{\epsilon}(q_RA)$
a continuous chain map which represents the corresponding bivariant
cohomology class in entire, analytic or local cyclic cohomology.
(One can choose a chain map which is independent of $R$ and natural in $A$). 
Let ${\cal E}'\,=\,({\cal H},\,\rho,\,\epsilon,\,F)$ be a bounded even
Fredholm module over $A$ and let 
$\psi:\,q_RA\to q_{C^*}A\to{\cal L(H)}$ (1.34) be its
characteristic homomorphism. Then the analytic cyclic cocycle 
$$
\check{ch}_{\cal P}({\cal
E}')\,=\,\iota^\epsilon_A\circ\psi_*\circ\tau_{\cal
P}\,\in CC^{anal}_0(A,\C)
\eqno(3.13)
$$
represents the abstract Chern-Connes character $\check{ch}$ (1.27) in
$K$-homology.

\item[b)] Let ${\cal E}(z)\,=\,({\cal
H},\,\rho,\,\epsilon,\,F(z)),\,z\in\C\backslash\R_-$ (1.5) be the holomorphic
family of bounded Fredholm modules associated to the given $\theta$-summable
unbounded Fredholm module ${\cal E}\,=\,({\cal H},\,\rho,\,\epsilon,\,{\cal
D})$ over $A$. Then the cochain 
$$
Ch_{\cal P}({\cal
E})\,=\,\frac{1}{\sqrt{\pi }\,i}\int_{1-i\infty}^{1+i\infty}\,\check{ch}_{\cal
P}({\cal E}(z^2))\, e^{z^2}\,dz
\eqno(3.14)
$$
extends to an entire cyclic cocycle over $A$.

\item[c)] The entire cyclic cocycle $Ch_{\cal P}({\cal E})$ represents 
the cohomology class of Connes' character $Ch_\epsilon({\cal E})$ 1.5 in 
the entire cyclic cohomology $HC_\epsilon^0(A)$ of $A$:
$$
[Ch_{\cal P}({\cal
E})
]\,=\,[Ch_\epsilon({\cal E})]\,\in\,HC_\epsilon^0(A)
\eqno(3.15)
$$
\end{itemize}
\end{theorem}

\begin{proof}

a)  The characteristic homomorphism (1.34) $\psi:\,qA\to{\cal K(H)}$ extends
to a continuous homomorphism $\psi:\,q_RA\to{\cal K(H)}$ for $R>\parallel
F\parallel$. The pullback of the analytic cyclic cocycle  $\vartheta$ along
$\psi$ and $\iota^\epsilon_A$ yields therefore an analytic cyclic
cocycle on $A$. As $\tau_{{\cal P}}$ represents the canonical generator of
$HC_{loc}^0({\cal K(H)})$ in local cyclic cohomology by (1.36) the formula
(1.39)  implies that the analytic cocycle $\check{ch}_{\cal P}$ represents the
Chern character in $K$-homology.

b) By (3.14) and the Cauchy integral formula (see \cite{Co}, pp.541-542)
$$
Ch_{\cal P}({\cal
E})\,=\,\frac{1}{\sqrt{\pi }\,i}\int_{1-i\infty}^{1+i\infty}\,\check{ch}_{\cal
P}({\cal E}(z^2))\, e^{z^2}\, dz\,=
$$
$$
=\,\iota^\epsilon_A\circ\,\frac{1}{\sqrt{\pi }\,i}\int_{1-i\infty}^{1+i\infty}
\psi_*\circ\tau_{{\cal P}}(z^2)\,d\mu(z)\,
=\,\iota^\epsilon_A\circ\beta
$$
and the latter cocycle is entire by 3.7.

c) Let $\pi:{\cal T}A\to A$ be the canonical epimorphism. Then for 
$R>>0$ sufficiently large
$$
\pi_*\circ Ch_{\cal P}({\cal
E})\,=\,\pi_*\circ\iota^\epsilon_A\circ\beta\,=
\,\iota^\epsilon_{{\cal T}A}\circ q_R(\pi)_*\circ\beta\,=\,
\iota^\epsilon_{{\cal T}A}\circ\nu_0 
$$
by 3.7. According to 3.5 and 3.6 the latter cocycle is cohomologous to
$\iota^\epsilon_{{\cal T}A}\circ\nu_1\,=\,\iota^\epsilon_{{\cal
T}A}\circ\widetilde{\tau}_{{\cal P}}$ which finally is cohomologous to
$\pi_*\circ Ch_{\epsilon}({\cal E})$ by 2.29. Therefore 
$$
\pi_*\circ [Ch_{\cal P}({\cal E})]\,=\,\pi_*\circ [Ch_{\epsilon}({\cal
E})]\,\in\, HC^\epsilon_0({\cal T}A,\C)
$$
The topological version of Goodwillie's theorem \cite{Pu2}, 5.1 asserts that 
the element $\pi_*\in HC^\epsilon_0({\cal T}A, A)$ is a bivariant
entire cyclic cohomology equivalence. So we deduce from the previous
equation
$$
[Ch_{\cal P}({\cal E})]\,=\,[Ch_{\epsilon}({\cal E})]\,\in\,
HC_\epsilon^0(A,\C)
$$
which was our claim.
\end{proof}

\section{Comparison of characters}

\begin{theorem}
Let ${\cal E}\,=\,({\cal H},\,\rho,\,\epsilon,\,{\cal D})$ be an even unbounded
$\theta$-summable Fredholm module over $A$. Let $\overline{A}$
be the enveloping $C^*$-algebra of $A$ and denote by \\ $[{\cal
E}]\,\in\,K^0(\overline{A})$ be the $K$-homology class of $\cal E$. 
Let $[Ch_\epsilon({\cal E})]\in
HC_\epsilon^0(A)$  be the cohomology class of
Connes' character cocycle 1.5 of $\cal E$ and let
$\check{ch}([{\cal E}])\,\in\,HC_{loc}^0(\overline{A})$ be the abstract 
Chern-Connes character (1.27) of the $K$-homology class of $\cal E$. Let
finally  $HC_\epsilon\,\to\,HC_{loc}$ be the canonical natural transformation
(\cite{Pu2}, section 7) from entire to local cyclic (co)homology. Then the
images of the characters  
$$ 
\begin{array}{ccccc}
HC_\epsilon^0(A) & \to & HC_{loc}^0(A) & \leftarrow &
HC_{loc}^0(\overline{A}) \\
  & & & & \\
Ch_{\epsilon}({\cal E}) & \to &  * & \leftarrow & \check{ch}([{\cal E}]) \\
\end{array} 
\eqno(4.1)
$$
in $HC_{loc}^0(A)$ coincide.
\end{theorem}

\begin{proof}
According to the previous theorem (see the proof of part b))
$$
Ch_{{\cal P}}({\cal E})\,=\,\frac{1}{\sqrt{\pi}i}\underset{\Gamma}{\int}
\check{ch}_{\cal P}(F(z^2))\, d\mu(z)
$$
Whereas only the whole integral defines an entire cyclic cocycle the
individual  terms $\check{ch}_{\cal P}(F(z^2))
\,=\,\iota^\epsilon_{A}\circ\psi_*\circ\tau_{\cal P}(z^2)$ make
sense as analytic cyclic cocycles. Variation  of the parameter $z\in\Gamma$ to
$z_0=1$ along $\Gamma$ gives rise to a smooth homotopy of the characteristic
homomorphisms (1.32). The Cartan homotopy formula in analytic cyclic cohomology
\cite{Pu}, 5.21 shows then that  
$$ 
\check{ch}_{\cal P}(F(z^2))\,=\,\check{ch}_{\cal
P}(F(1))\,+\,\partial(\chi(z^2)) 
\eqno(4.2)
$$
where $\chi(z^2)\in CC^1_{anal}(A)$ is an analytic cyclic cochain over
$A$. The growth conditions of the cochain $\chi(z^2)$ are satisfied uniformly 
 with respect to $z\in \Gamma$ as a look at 2.26, 2.27 and the definition of
$\chi$ shows. Therefore 
$$
\chi\,=\,\frac{1}{\sqrt{\pi}i}\underset{\Gamma}{\int}\chi(z^2)\,d\mu(z)
$$ 
extends to an analytic cyclic cochain on $A$. Integration of (4.2) over
$\Gamma$ provides the equality  
$$
Ch_{{\cal P}}({\cal E})\,=\,\check{ch}_{\cal
P}(F(1))\,\frac{1}{\sqrt{\pi}i}\underset{\Gamma}{\int}d\mu(z)\,+\,\partial(\chi)\,=\,
\check{ch}_{\cal P}(F(1))\,+\,\partial(\chi)
$$
of analytic cyclic cocycles.
If we pass thus from entire to analytic cyclic cohomology \cite{Me},\cite{Pu}
 we find that 
$$
[Ch_{\epsilon}({\cal E})]\,=\,[Ch_{{\cal P}}({\cal E})]\,=\,[\check{ch}_{\cal
P}(F(1))]\in HC_{anal}^0(A) 
$$
coincide. Let $\overline{A}$ be the enveloping $C^*$-algebra of $A$ and
consider the commutative diagram (in the notations of (1.33) and (2.43))
$$
\begin{array}{ccccccccc}
0 & \to & {\cal K(H)} & \longrightarrow &  {\cal L(H)} & \longrightarrow
& {\cal Q(H)} & \to & 0 \\
 & & \psi\uparrow & & \uparrow & & \uparrow & & \\
0 & \to & q_{C^*}\overline{A} & \longrightarrow & Q_{C^*}\overline{A} &
\longrightarrow & \overline{A} & \to & 0 \\
 & & \uparrow & & \uparrow & & \uparrow i & & \\
0 & \to & q_RA & \longrightarrow & Q_RA & \longrightarrow
& A & \to & 0 \\
\end{array}
$$
which exists for all $R\geq 1$ as homomorphisms of $C^*$-algebras are
contractive. A look at the formula (3.14) for the Chern-Connes character
shows that 
$\check{ch}_{\cal P}(F(1))\in CC_{anal}^0(A)$ extends to an analytic
cyclic cocycle on the $C^*$-algebra $\overline{A}$. Viewed as a local cyclic
cocycle it represents the Chern-Connes character 
$\check{ch}({\cal E}(1))\in HC^0_{loc}(\overline{A})$ of the $K$-homology class
of the bounded Fredholm module ${\cal E}(1)=({\cal H},\rho,\epsilon, F(1))$.
Now by definition  
$$
F(1)\,=\,\frac{{\cal D}+\epsilon}{(1+{\cal D}^2)^{\frac12}}\,=\,f({\cal
D})\,+\,\epsilon\,(1+{\cal D}^2)^{-\frac12},\,\,\,f(x)\,=\,x(1+x^2)^{-\frac12}
$$
Therefore $f({\cal D})$ is a compact perturbation of $F(1)$ which shows that
this bounded Fredholm module represents the $K$-homology class
of ${\cal E}$. So we deduce finally 
$$
\check{ch}([{\cal E}])\,=\,[\check{ch}_{\cal
P}(F(1))]\,=\,[Ch_{{\cal P}}({\cal
E})]\,=\,[Ch_{\epsilon}({\cal E})]\,
\in\,HC^0_{loc}(\overline{A})  
$$
which proves the theorem.
\end{proof}

\section{On a problem of Connes}

As an application of the comparison result for characters of Fredholm modules 
we present a partial solution of a problem posed by A.~Connes in \cite{Co3},
p.83 (see also \cite{Co2}, pp.407-414). There he constructs a certain
infinite dimensional unbounded $\theta$-summable Fredholm
module over the reduced group $C^*$-algebra of a lattice in a semisimple
Lie-group. The $K$-homology class of this module is closely related to
Kasparov's $\gamma$-element \cite{Ka1}, \cite{Co2}. Connes predicts that the
character cocycle of this Fredholm module is cohomologous to the canonical
trace and asks to verify this. We will determine the class of
Connes' character cocycle in local cyclic cohomology for uniform lattices in
semisimple Lie groups of real rank one. This is as far as I know the
first calculation of the Connes character of an infinite dimensional unbounded
$\theta$-summable Fredholm module.

We begin by a presentation of Connes' problem.
\begin{definition}\cite{Co3},pp.81-82, \cite{Co2},pp.408-409.

Let $\Gamma$ be a discrete subgroup of a real semisimple Lie group $G$. Let 
$G/K$ be the associated symmetric space and let $d(-,-)$ be the distance
function associated to a (left)-translation invariant Riemannian metric on
$G/K$. Let $L^2(G/K,\,\Lambda^*)$ be the Hilbert space of square-integrable 
differential forms graded by the decomposition into forms of even and odd
degree, respectively. Fix a base point $x_0\in G/K$. Then 
$$
\begin{array}{cc}
{\cal D}_{(G,x_0)}\,=\,e^\varphi\cdot d\cdot e^{-\varphi}\,+\,(e^\varphi\cdot
d\cdot e^{-\varphi})^*, & \varphi(x)\,=\,d(x,x_0)^2 \\
\end{array}
\eqno(5.1)
$$
is a densely defined, unbounded, odd, essentially selfadjoint operator on
$L^2(G/K,\,\Lambda^*)$ and 
$$
{\cal E}_\Gamma\,=\,\left(L^2(G/K,\,\Lambda^*),\,{\cal D}_{(G,x_0)}\right)
\eqno(5.2)
$$
defines a $\theta$-summable unbounded even Fredholm module over a dense 
and holomorphically closed Banach subalgebra ${\cal A}(\Gamma)\subset
C^*_r(\Gamma)$ of the reduced group $C^*$-algebra $C^*_r(\Gamma)$ of $\Gamma$.
\end{definition}

This Fredholm module cannot be finitely summable according to \cite{Co5}
because $G$ contains non-amenable discrete subgroups.

The problem which Connes poses in \cite{Co3}, p.83 is equivalent to the

\begin{problem}{\bf (Connes)}

Show that the Connes character of ${\cal E}_\Gamma$ is cohomologous to the restriction 
of the canonical trace $\tau_{can}$ on $C^*_r(\Gamma)$ to ${\cal A}(\Gamma)$:
$$
[Ch_\epsilon({\cal E}_\Gamma)]\,=\,[\tau_{can}]\,\in\,HC_\epsilon^0({\cal
A}(\Gamma)) \eqno(5.3)
$$
\end{problem}

The interest in this problem stems from the fact that a positive answer
provides a solution of the Kadison-Kaplansky idempotent conjecture for the
considered group. 

We recall what is known about the $K$-homology class of
${\cal E}_\Gamma$. The next lemma collects material from \cite{Ka1} and
\cite{Co2}.

\begin{lemma}
Let ${\cal E}_\Gamma$ be the unbounded Fredholm module introduced in 5.1 and
let $[{\cal E}_\Gamma]\,\in\,K^0(C^*_r(\Gamma))$ be its associated $K$-homology
class. Denote by $p:C^*_{max}(\Gamma)\to C^*_r(\Gamma)$ the canonical
epimorphism. Let 
$$
KK(-\rtimes_{max}\Gamma,\C)\,\simeq\,KK^\Gamma(-,\C)
\eqno(5.4)
$$ 
be the tautological isomorphism. In particular
$$KK(C^*_{max}\Gamma,\C)\,\simeq\,KK^\Gamma(\C,\C)$$
Then the pullback $p_*\circ[{\cal
E}_\Gamma]\,\in\,KK_0(C^*_{max}(\Gamma),\C)$ 
corresponds under  the tautological isomorphism 
to Kasparov's canonical $\gamma$-element
$\gamma_\Gamma\,\in\,KK^\Gamma(\C,\C)$
$$
\begin{array}{ccccc}
KK(C^*_{r}(\Gamma),\C) & \longrightarrow & KK(C^*_{max}(\Gamma),\C) & \simeq &
KK^\Gamma(\C,\C) \\
 & & & & \\
{[}{\cal E}_\Gamma{]} & \longrightarrow  & p_*\circ{[}{\cal
E}_\Gamma{]} & \leftrightarrow & \gamma_\Gamma \\
\end{array}
\eqno(5.5)
$$
\end{lemma}

The abstract characterization of the considered Fredholm module 
will allow us to apply, according to the comparison result of the previous
section, our partial calculation of the abstract Chern-Connes character 
of the $\gamma$-element in \cite{Pu3}. This will lead to the promised partial
solution of Connes' problem.

In the next step we clarify the behaviour of the tautological isomorphism
(5.4)  with respect to the various bivariant Chern-Connes characters.

\begin{lemma}
Let $\Gamma$ be a countable discrete group, let $A$ be a separable
$\Gamma$-$C^*$-algebra and let
$\,A\rtimes_{max}\Gamma$ be the associated maximal crossed
product $C^*$-algebra. Let
$\iota:\,\ell^1(\Gamma,A)\to A\rtimes_{max}\Gamma$ be the canonical inclusion
and denote by $\pi:\ell^1(\Gamma)\,\to\,\C$ the trivial representation.
Let
$$
ch_{biv}:\,KK_*(-,-)\,\longrightarrow\,HC_*^{loc}(-,-)
\eqno(5.6)
$$
be the bivariant Chern-Connes character \cite{Pu1}, section 6 and let 
$$
ch^\Gamma_{\ell^1(\Gamma)}:\,KK^\Gamma(-,-)\,\longrightarrow\,
HC_0^{loc}(\ell^1(\Gamma,-),\ell^1(\Gamma,-))
\eqno(5.7)
$$
be the eqivariant bivariant Chern-Connes character \cite{Pu3}, 5.2.
Then there is a natural commutative diagram
$$
\begin{array}{ccc}
 KK^\Gamma(A,\C)  & \simeq &
 KK(A\rtimes_{max}\Gamma,\C)\\  
 & & \\  
 & &  \downarrow ch_{biv} \\
 & & \\
 ch^\Gamma_{\ell^1(\Gamma)} \downarrow & & HC_0^{loc}(A\rtimes_{max}\Gamma,\C)\\  
 & & \\
  & &  \downarrow \iota_*\circ -\\
 & & \\
HC_0^{loc}(\ell^1(\Gamma,A),\ell^1(\Gamma)) &
\overset{-\circ\pi_*}{\longrightarrow} & HC_0^{loc}(\ell^1(\Gamma,A),\C)
\\  
\end{array}
\eqno(5.8)
$$
\end{lemma}

\begin{proof}

Observe that the diagram is natural in the following sense. For every
$\Gamma$-$C^*$-algebra $B$ the abelian group $KK^\Gamma(B,A)$ acts on 
the groups in the diagram by left multiplication, and all arrows in the
diagram are compatible with this action. In fact, any
element $\alpha\in KK^\Gamma(B,A)$ acts via left multiplication by
$j(\alpha)\in KK(B\rtimes_{max}\Gamma,A\rtimes_{max}\Gamma)$  on
$KK(A\rtimes_{max}\Gamma,-)$, via $ch_{biv}(j(\alpha))$ on
$HC_0^{loc}(A\rtimes_{max}\Gamma,-)$, and via
$ch^\Gamma_{\ell^1(\Gamma)}(\alpha)$ on $HC_0^{loc}(\ell^1(\Gamma,A),-)$.  The
compatibility of these actions with the arrows in the diagram is equivalent  to
the multiplicativity of the descent homomorphism $j$ \cite{Ka1} and 
the multiplicativity of the various bivariant Chern-Connes characters
\cite{Pu1}, \cite{Pu3}, respectively. Every element $\alpha\in
KK^\Gamma(A,\C)$ can be written as a product $\alpha\,=\,\alpha\otimes
1_\Gamma$  where $1_\Gamma\in KK^\Gamma(\C,\C)$ is the unit. Therefore it
suffices to verify that the diagram commutes
when applied to the unit class $1_\Gamma\in KK^\Gamma(\C,\C)$. Now the unit
class corresponds  under the tautological isomorphism (5.4) to the trivial
representation of $C^*_{max}(\Gamma)$. The assertion follows then from the
naturality of the Chern-Connes characters.

\end{proof}

The main result of \cite{Pu3} yields by the previous lemma 

\begin{theorem}
Let $\Gamma$ be a torsion free, cocompact, discrete subgroup of 
a semisimple Lie group $G$ of real rank one. Let 
$$
[{\cal E}_\Gamma]\,\in\,K^0(C^*_r(\Gamma))\,=\,KK(C^*_r(\Gamma),\C)
$$
be the $K$-homology class of Connes' Fredholm module ${\cal E}_\Gamma$ 
introduced in 5.1. Then its abstract Chern-Connes character is given by the
cohomology class  of the canonical trace 
$$
\check{ch}([{\cal E}_\Gamma])\,=\,[\tau_{can}]\,\in\,HC^0_{loc}(C^*_r(\Gamma))
\eqno(5.9)
$$
\end{theorem}

\begin{proof}
We proceed in several steps. 

1) According to our assumptions the space $\Gamma \backslash G/K$ is a compact
Riemannian manifold of negative sectional curvature. Therefore the
group $\Gamma$ is word-hyperbolic. In particular, it satisfies  Jolissaint's
condition of rapid decay \cite{CM},section 6. This means that the space
${\mathfrak A}(\Gamma)$ of square-summable functions of rapid decay (w.r.t. a
word metric) on $\Gamma$ forms an algebra under convolution which is dense and
holomorphically closed in $C^*_r(\Gamma)$. 

2) It was shown in \cite{Pu3} that for a word-hyperbolic group $\Gamma$ the
local cyclic cohomology groups of the algebras $\ell^1(\Gamma),\,{\mathfrak
A}(\Gamma)$  and of $\ell^1(\Gamma)\cap{\mathfrak A}(\Gamma)$ decompose as the
direct  sum of a homogeneous and an inhomogeneous part. The homogeneous
part corresponds to the contribution of the conjugacy class of the
unit and the inhomogeneous part to the contribution of the other conjugacy
classes of $\Gamma$. Moreover the various inclusions of these algebras induce
isomorphisms
$$
HC^*_{loc}({\mathfrak A}(\Gamma))_{hom}\,\overset{\simeq}{\longrightarrow}
HC^*_{loc}(\ell^1(\Gamma)\cap{\mathfrak
A}(\Gamma))_{hom}\,\overset{\simeq}{\longleftarrow} 
HC^*_{loc}(\ell^1(\Gamma))_{hom} 
\eqno(5.10)
$$
of the homogeneous parts of the corresponding local cyclic cohomology groups.

3) Recall that by definition the abstract
character $\check{ch}$ in $K$-homology corresponds to the bivariant
Chern-Connes character $ch_{biv}$ under the identification
$K^*(-)\,=\,KK_*(-,\C)$. We calculate the pull back 
of $\check{ch}([{\cal E}_\Gamma])\,\in\,HC^0_{loc}(C^*_r(\Gamma))$  along the
inclusion  $i:\,\ell^1(\Gamma)\,\to\,C^*_r(\Gamma)$. Lemma 5.3 and the
commutative diagram of lemma 5.4 show that  
$$  
i_*\circ\check{ch}([{\cal
E}_\Gamma])\,=\,\iota_*\circ p_*\circ\check{ch}([{\cal
E}_\Gamma])\,=\,\iota_*\circ\check{ch}(p^*[{\cal
E}_\Gamma])\,=\,ch^\Gamma_{\ell^1(\Gamma)}(\gamma_\Gamma)\circ\pi
\eqno(5.11)
$$
In \cite{Pu3} we gave a qualitative characterization of the equivariant
Chern-Connes character of the
$\gamma$-element for a torsion free word-hyperbolic group:
$$
ch^\Gamma_{\ell^1(\Gamma)}(\gamma_\Gamma)\,=\,\pi_{Hom}\in 
\,HC_0^{loc}(\ell^1(\Gamma),\ell^1(\Gamma))
\eqno(5.12)
$$ 
where $\pi_{Hom}$ denotes the canonical projection onto the homogeneous part 
of the local cyclic complex of $\ell^1(\Gamma)$. So we conclude
$$
i_*\circ\check{ch}([{\cal
E}_\Gamma])\,=\,\pi_{Hom}\circ\pi_*\,=\,i_*\circ[\tau_{can}]\in 
\,HC^0_{loc}(\ell^1(\Gamma))
\eqno(5.13)
$$

4) Let $i':{\mathfrak A}(\Gamma)\,\to\,C^*_r(\Gamma)$ be the canonical
inclusion. By 2) we have 
$$
i'_*\circ\check{ch}([{\cal
E}_\Gamma])\,=\,(i'_*\circ\check{ch}([{\cal
E}_\Gamma]))_{hom}\,+\,\,(i'_*\circ\check{ch}([{\cal
E}_\Gamma]))_{inhom}\,\in\,HC^0_{loc}({\mathfrak A}(\Gamma))
$$
and from 2) and 3) we learn
$$
(i'_*\circ\check{ch}([{\cal
E}_\Gamma]))_{hom}\,=\,i'_*\circ[\tau_{can}]\,\in\,HC^0_{loc}({\mathfrak
A}(\Gamma))_{hom}  
\eqno(5.14)
$$

5) We calculate the inhomogeneous part of $i'_*\circ\check{ch}([{\cal
E}_\Gamma])$. Recall that the $\gamma$-element appears as a Kasparov product 
$$
\gamma\,=\beta\otimes\alpha,\,\beta\in KK^\Gamma(\C,C_0(X)),\,\alpha\in
KK^\Gamma(C_0(X),\C) 
\eqno(5.15)
$$
where $X$ is a manifold on which $\Gamma$ acts properly and freely. Because the
action of $\Gamma$ on $X$ is proper, the associated maximal and reduced crossed
products coincide:
$$
C_0(X)\rtimes_{max}\Gamma\,\overset{\simeq}{\longrightarrow}\,C_0(X)\rtimes_r\Gamma
\eqno(5.16)
$$  
The compatibility of the tautological isomorphism (5.4) with the Kasparov
product, discussed in the proof of 5.4, shows that
$$
[{\cal E}_\Gamma]\,=\,j_r(\beta)\otimes\alpha',\,j_r(\beta)\in
KK(C^*_r\Gamma,C_0(X)\rtimes_r\Gamma),\,\alpha'\in KK(C_0(X)\rtimes_r\Gamma,\C)
$$
where the class $\alpha'\in
KK(C_0(X)\rtimes_r\Gamma,\C)=KK(C_0(X)\rtimes_{max}\Gamma,\C)$ corresponds to 
$\alpha\in KK^\Gamma(C_0(X),\C)$ under (5.4). The multiplicativity of the
Chern-Connes character implies
$$
i'_*\circ\check{ch}([{\cal
E}_\Gamma])\,=\,i'_*\circ ch_{biv}(j_r(\beta))\circ\check{ch}(\alpha') 
\,=\,ch^\Gamma_{{\mathfrak A}(\Gamma)}(\beta)
\circ i'_*\circ\check{ch}(\alpha')
\eqno(5.17)
$$
$ch^\Gamma_{{\mathfrak A}(\Gamma)}(\beta)\in HC_0^{loc}({\mathfrak
A}(\Gamma),{\mathfrak A}(\Gamma,C_0(X)))$, $i'_*\circ\check{ch}(\alpha')\in
 HC^0_{loc}({\mathfrak A}(\Gamma,C_0(X)))$, and
similar identities hold for the homogeneous and inhomogeneous parts
\cite{Pu3}, 5.2, respectively. It has been shown in \cite{Pu3}, 4.8 that  
$$
HC^0_{loc}({\mathfrak A}(\Gamma,C_0(X)))_{inhom}\,=\,0
\eqno(5.18)
$$ 
because $\Gamma$ acts
properly on $X$. We deduce from this and the previous identity that
$$
\begin{array}{ccc}
i'_*\circ\check{ch}(\alpha')_{inhom}=0 & \text{and} & 
i'_*\circ\check{ch}([{\cal
E}_\Gamma])_{inhom}\,=\,0
\end{array}
\eqno(5.19)
$$

6) From parts 4) and 5) we may conclude that
$$
i'_*\circ\check{ch}([{\cal
E}_\Gamma])\,=\,i'_*\circ\tau_{can}\,\in\,HC^0_{loc}({\mathfrak A}(\Gamma))
\eqno(5.20)
$$
Cowling and Haagerup have shown \cite{CH}, that the reduced group
$C^*$-algebra $C^*_r(\Gamma)$ of a cocompact discrete subgroup of a semisimple
Lie group of real rank one possesses the Grothendieck approximation property.
This implies that the inclusion of any dense subalgebra ${\cal
R}\hookrightarrow C^*_r(\Gamma)$, which is the domain of an unbounded
derivation on $C^*_r(\Gamma)$, induces a local cyclic cohomology equivalence
\cite{Pu3}, 5.15. and 3.10. The Jolissaint algebra is of this type so that we
deduce that  $$
i'_*\,\in\,HC_0^{loc}({\mathfrak A}(\Gamma),C^*_r(\Gamma))
\eqno(5.21)
$$
is a local cyclic cohomology equivalence for the groups under consideration. 
Thus the identity (5.20) implies 
$$
\check{ch}([{\cal
E}_\Gamma])\,=\,[\tau_{can}]\,\in\,HC^0_{loc}(C^*_r(\Gamma))
\eqno(5.22)
$$
which was to be proved.
\end{proof}

As a consequence we obtain finally a partial solution of Connes' problem.

\begin{theorem}
Let $\Gamma$ be a torsion free, cocompact, discrete subgroup of 
a semisimple Lie group $G$ of real rank one. Let 
$$
{\cal E}_\Gamma\,=\,\left(L^2(G/K,\,\Lambda^*),\,{\cal D}_{(G,x_0)}\right)
$$
be Connes' unbounded $\theta$-summable Fredholm module over the dense 
and holomorphically closed Banach subalgebra ${\cal A}(\Gamma)\subset
C^*_r(\Gamma)$.
Then Connes' character cocycle 
$$
Ch_\epsilon({\cal E}_\Gamma)\,\in\,CC_\epsilon^0({\cal A}(\Gamma))
$$
is cohomologous in local cyclic cohomology to the canonical trace 
on ${\cal A}(\Gamma)$:   
$$
[Ch_\epsilon({\cal E}_\Gamma)]\,=\,[\tau_{can}]\,\in\,HC^0_{loc}({\cal
A}(\Gamma)) 
\eqno(5.23)
$$
\end{theorem}

\begin{proof}
Let $j:\,{\cal A}(\Gamma)\,\to\,C^*_r(\Gamma)$ be the inclusion. According
to the comparison theorem 4.1 and theorem 5.5 the image of Connes' character
cocycle $Ch_\epsilon({\cal E}_\Gamma)$ in local cyclic cohomology satisfies
$$
[Ch_\epsilon({\cal E}_\Gamma)]\,=\,j_*\circ\check{ch}([{\cal
E}_\Gamma])\,=\,j_*\circ[\tau_{can}]\,=\,[\tau_{can}]\,\in\,HC^0_{loc}({\cal
A}(\Gamma))
\eqno(5.24) 
$$
\end{proof}

As a well known consequence (see \cite{Co1}, \cite{Co2}) we note the following
result which was obtained  previously with other methods by Lafforgue
\cite{La}, Mineyev-Yu \cite{MY}, Kasparov-Yu (unpublished), and myself
\cite{Pu3}.

\begin{theorem} (Kadison-Kaplansky conjecture)

Let $\Gamma$ be a torsion free, cocompact, discrete subgroup of 
a semisimple Lie group $G$ of real rank one. Then the reduced group 
$C^*$-algebra $C^*_r(\Gamma)$ contains no idempotents except 0 and 1.
\end{theorem}

\begin{proof} (following Connes)
By a well known argument it suffices to show that the canonical trace 
takes only integer values on idempotents in $C^*_r(\Gamma)$. Let 
$e=e^2\in {\cal A}(\Gamma)$ and let $[e]\in K_0({\cal A}(\Gamma))$ be its class
in $K$-theory. The index formula 1.6 and the previous theorem show that
$$
\tau_{can}(e)\,=\,\langle [\tau_{can}], ch([e])\rangle\,=\,
\langle [Ch_\epsilon({\cal E}_\Gamma)], ch([e])\rangle\,=\,Ind({\cal
E}_\Gamma\otimes [e])\,\in\,\Z
\eqno(5.25)
$$
Because ${\cal A}(\Gamma)$ is dense and holomorphically closed in
$C^*_r(\Gamma)$ the same conclusion holds for idempotents $e=e^2\in
C^*_r(\Gamma)$. This implies the theorem. 
\end{proof}

\end{document}